
%

%
%

%
%

\font\eusm=eusm10                          

\font\eusms=eusm7                          

\font\eusmss=eusm5                         

\font\scriptsize=cmr7

\font\tiny=cmr5

\input amstex

\documentstyle{amsppt}
  \magnification=1000
  \hsize=7.0truein
  \vsize=9.0truein
  \hoffset -0.1truein
  \voffset 0.75truein
  \parindent=2em

\define\Ac{{\Cal A}}                           

\define\brakE#1{                               
     \left[ #1 \right]_{E}}

\define\brakphi#1{                             
     \left[ #1 \right]_{\phi}}

\define\brakr#1{                               
     \left[ #1 \right]_{r}}

\define\brakstar#1{                            
     \left[ #1 \right]_{\ast}}

\define\Case#1#2{\newline\noindent             
     {\bf Case #1: #2.} \newline}

\define\Cc{\text{\eusm C}}                     

\define\Cpx{\bold C}                           

\define\dif{\text{\it d}}                      

\define\Ec{\text{\eusm E}}                     

\define\freeF{\bold F}                         

\define\freeprod{                              
     \operatornamewithlimits{\ast}}

\define\ih{\hat{\imath}}                       

\define\Imaginary{\text{\rm Im}\,}             

\define\Integers{\bold Z}                      

\define\iotat{\tilde{\iota}}                   

\define\jh{\hat{\jmath}}                       

\define\MnC{M_n(\Cpx)}                         

\define\MNC{M_N(\Cpx)}                         

\define\MvN{{\Cal M}}                          

\define\Naturals{\bold N}                      

\define\ngi{{n\rightarrow\infty}}              

\define\nbyn{n\times n}                        

\define\OdO{0\cdots0}                          

\define\Orient{\text{\eusm O}}                 

\define\Pc{\text{\eusm P}}                     

\define\QED{$\hfill$\qed\enddemo}              

\define\Qg{\mathchoice                         
     {\text{\eusm Q}}
     {\text{\eusm Q}}
     {\text{\eusms Q}}
     {\text{\eusmss Q}} }

\define\Rlpt{\text{\rm Re}\;}                  

\define\Reals{\bold R}                         

\define\restrict{\lower .3ex                   
     \hbox{\text{$|$}}}

\define\Rl{\text{\eusm R}}                     

\define\Rtilde{\overset                        
     \text{\eusms R}\to\sim}

\define\Step#1{\noindent                       
     {\bf Step #1.} $\,$ }

\define\vtilde{\overset v\to\sim}              


\topmatter

  \title On certain free product factors
           via an extended matrix model \endtitle

  \author Ken Dykema \endauthor

  \affil University of California, \\
         Berkeley, California, USA 94720, \\
         (e--mail dykema\@math.berkeley.edu) \endaffil

  \thanks{Studies and research
          supported by the Fannie and John Hertz Foundation;
          this work supported in part by the Natural Sciences and
          Engineering Council of Canada and by the funds
          FCAR du Qu\'ebec.} \endthanks
  \thanks{This work will form part of the author's Ph.D\. thesis
          at the University of California at Berkeley} \endthanks
  \date{ November 1991 } \enddate

  \abstract{
    A random matrix model for freeness is extended and used to investigate
      free products of free group factors with matrix algebras and with
      the hyperfinite II$_1$--factor.
    The latter is shown to be isomorphic to a free group factor having one
      additional generator.} \endabstract

\endtopmatter

\document
  \TagsOnRight
  \baselineskip=14pt

\noindent{\bf Introduction.}

  The finite von Neumann algebra $L(G)$ associated to a discrete
    group $G$ was introduced by Murray and von Neumann~\cite{5}.
  It is the von Neumann algebra generated by the representation
    of $G$ on $l^2(G)$ by left translation operators, with faithful
    trace given by the vector--state for $\delta_e$.
  They gave the free group factors $L(\freeF_K)$
    (where $\freeF_K$ is the nonabelian free group on $K$ generators),
    ($2\le K\le\infty$) as examples of type~II$_1$ factors which are
    not hyperfinite.
  Much of the structure around free group factors is not yet
    understood.
  For example, the isomorphism question of
    whether $L(\freeF_n)\cong L(\freeF_m)$ for $n\neq m$ remains open.
  In another direction, one can ask:
   for which groups $G$ is $L(G)$ a free group factor?

  Voiculescu \cite{8,9,10,11,12} has created a theory of freeness
    in noncommutative probability spaces, which has begun to
    shed light around the free group factors~\cite{10,6}.
  In this theory there is a notion of the free product
    of finite von Neumann algebras with specified traces, (see also~\cite{2}),
    for which one has $L(G_1)*L(G_2)\cong L(G_1*G_2)$.

  Murray and von Neumann showed also in~\cite{5} that there is
    up to isomorphism
    only one hyperfinite factor of type~II$_1$,
    which we denote by $R$.
  In his fundamental paper~\cite{3},
    Alain Connes has shown that $L(G)\cong R$ if and only if $G$
    is an amenable i.c.c\. group.

  In this paper, we shall show that
    \roster
      \item"(a)" $L(\freeF_K)\ast M_N(\Cpx)
        \cong L(\freeF_{N^2K})\otimes M_N(\Cpx)
        \qquad(1\le K\le\infty,\,N\ge2)$
      \item"(b)" $L(\freeF_K)\ast R\cong L(\freeF_{K+1})
        \qquad(1\le K\le\infty)$.
    \endroster
  (Of course we have $\freeF_1=\Integers$.)
  This together with Connes' results implies, for example,
    that $L(\Integers*G)\cong L(\freeF_2)$
    whenever $G$ is an amenable i.c.c\. group.

  In order to prove~(a) and~(b), we use random matrices.
  Wigner~\cite{14,15} showed that certain self--adjoint $\nbyn$ random
    matrices with independent entries have distributions tending
    to a semicircle law as $\ngi$.
  (See~\cite{1,4,7,13} for a sample of the literature
    about random matrices.)
  In~\cite{11}, Voiculescu showed that families of such random matrices having
    mutually independent Gaussian entries become, together with diagonal
    matrices of constant random variables, asymptotically free in the
    limit as $\ngi$.
  (He then used this matrix model in~\cite{10} to prove certain
    isomorphisms around free group factors.)
  We will extend this matrix model simultaneously in two respects.
  First, we may replace the Gaussian entries of the matrices in
    Voiculescu's theorem with more general random variables,
    including all cases which Wigner considered.
  Secondly, we may replace Voiculescu's constant diagonal matrices
    with constant block diagonal matrices, such that the block size
    remains constant (or grows slowly) as $\ngi$.

  This paper has four sections.
  In~\S1, we give some preliminaries about freeness;
   in~\S2, we prove the extended matrix model;
   in~\S3 we prove~(a) and find one more algebra
    isomorphic to these two;
   in~\S4 we prove~(b).

\noindent{\bf \S1\. Preliminaries}

  We give here some aspects of Voiculescu's theory of freeness
    which we will use later.
  For more details, see~\cite{10,12} and other papers of
    Voiculescu.

  A {\it $\ast$--noncommutative probability space}\/
    is $(A,\phi)$, where $A$ is a
    unital $\ast$--algebra over~$\Cpx$ and $\phi$ a positive state on $A$
    sending $1$ to $1$.
  A W$^\ast$--noncommutative probability space is where $A$
    is a von Neumann algebra and $\phi$ is
    ultra--weakly continuous.
  (In addition, we will always have that $\phi$ is a faithful trace.)
  {\it Random variables}\/ are elements $a\in A$.
  The {\it distribution}\/ of $a\in A$ is the linear functional,
    $\mu_a:\Cpx[X]\rightarrow\Cpx$ which sends the polynomial, $P$,
    to $\phi(P(a))$.
  In the W$^\ast$ setting, if $a$ is normal,
    then $\mu_a$ is given by
    $$ \mu_a(P)=\int_{\sigma(a)}P(w)\dif\lambda(w), $$
    where $\dif\lambda$ is equal to $\phi$ of spectral measure, supported on
    the spectrum of $a$.
  A family, $(a_i)_{i\in I}$ of random variables in $A$ has a
    {\it joint distribution}\/, the linear functional
    $\mu:{\Cpx\langle X_i\mid i\in I\rangle}\rightarrow\Cpx$
    (where ${\Cpx\langle X_i\mid i\in I\rangle}$
    is the algebra of noncommutative
    polynomials), which is $\mu=\phi\circ h$, where
    $h:{\Cpx\langle X_i\mid i\in I\rangle}\rightarrow A$ is the unique unital
    homomorphism sending $X_i$ to $a_i$.
  A sequence, $\mu_n$, of (joint) distributions is said to converge
    to the {\it limit distribution,}\/ $\mu$, if
    $\lim_{n\rightarrow\infty}\mu_n(P)=\mu(P)$
    $\forall\,P\in{\Cpx\langle X_i\mid i\in I\rangle}$.

  A family, $(A_i)_{i\in I}$, of unital subalgebras of $A$
    is {\it free}\/ if $\phi(a_1\cdots a_n)=0$ whenever $a_j\in A_{i_j}$,
    $i_1\neq i_2\neq\cdots\neq i_n$ and $\phi(a_j)=0\;\forall j$.
  A family of subsets of $A$ (or of random variables of $A$)
    is said to be free if the unital subalgebras they generate are
    free, and {\it $\ast$--free}\/ if the unital $\ast$--algebras
    they generate are free.

  \proclaim{Proposition 1.1} {\rm(\cite{12})}.
    \roster
    \item"(i)" If $(A_i)_{i\in I}$ is a free family of subalgebras in
      $(A,\phi)$ and $(B_j)_{j\in J}$ are the subalgebras generated by
      $\cup_{i\in I_j}A_i$, where $I=\coprod_{j\in J}I_j$ is a partition
      of $I$, then the family $(B_j)_{j\in J}$ is free.
    \item"(ii)" If $(A_i)_{i\in I}$ is a free family of subalgebras in
      $(A,\phi)$ and $(C_{ik})_{k\in K_i}$ is a free family of subalgebras in
      $(A_i,\phi|_{A_i})$, then the family of subalgebras
      $(C_{ik})_{(i,k)\in K}$ is free in $(A,\phi)$,
      where $K=\coprod_{i\in I}\{i\}\times K_i$.
    \endroster
  \endproclaim

  A sequence, $(a_i(n))_{i\in I}$, for $n\in\Naturals$, of families of
    random variables is said to be {\it asymptotically free}\/ if their
    distributions converge to a limit distribution, $\mu$, and the
    random variables $(X_i)_{i\in I}$ in
    $(\Cpx\langle X_i\mid i\in I\rangle,\mu)$ are free.

  If $(A,\phi)$ and $(B,\psi)$ are $\ast$--noncommutative probability spaces,
    then there is a state, denoted $\phi\ast\psi$ on the
    free product $\ast$--algebra, $A\ast B$, such that, $A$ and $B$ being
    identified with canonical subalgebras of $A\ast B$, we have
    $\phi\ast\psi|_A=\phi$ and $\phi\ast\psi|_B=\psi$, and
    $A$ and $B$ are free in $(A\ast B,\phi\ast\psi)$.
  If $A$ and $B$ are von Neumann algebras, the ultra--weak closure of
    the GNS--representation of $(A\ast B,\phi\ast\psi)$ is what we call
    the {\it free product}\/ of the von Neumann algebras (with specified
    states).

  Note that, since a von Neumann algebra, $M$, with faithful trace,
    $\phi$, can be faithfully represented on $L^2(M,\phi)$, its
    isomorphism class is completely specified by the joint
    distribution of a set of generators of $M$ (together with their adjoints).
  Thus, for example, (utilizing the functional calculus) we see that
    any such von Neumann algebra generated by a $\ast$--free set of
    $K$ normal elements whose traces of spectral measure
    have no atoms is $L(\freeF_K)$.

  A random variable, $a$, in a $\ast$--noncommutative probability space is
    said to be {\it semicircular}\/ (of radius $r$) if it is
    self--adjoint and its distribution is given by
    $$ \mu_a(X^m)=\frac{2}{\pi r^2}\int_{-r}^rt^m(r^2-t^2)^{1/2}\dif t
      \;\;\;\forall m\ge0. $$
   it is said to be {\it quarter--circular}\/ (of radius $r$)
    if it is positive and its
    distribution is given by
    $$ \mu_a(X^m)=
      \frac{4}{\pi r^2}\int_0^rt^m(r^2-t^2)^{1/2}\dif t\;\;\;\forall m\ge0 $$
   it is said to be {\it circular}\/ (of radius $r$) if
    $\left\{\frac{a+a^\ast}{\sqrt{2}},\,\frac{a-a^\ast}{i\sqrt{2}}\right\}$
    is a free pair of random variables, each semicircular of radius $r$.
  A {\it Haar unitary}\/ is a unitary $u$ in $A$ such that
    $\phi(u^k)=0$ $\forall k\neq0$.
  This is equivalent to $\phi$ of the spectral measure of $u$
    being Haar measure on the circle.

  \proclaim{Theorem 1.2}{\rm (\cite{10})}.
    In a $\ast$--noncommutative probability space $(A,\phi)$, where $A$ is
      a von Neumann algebra and $\phi$ a faithful trace, $a\in A$
      is circular if and only if its polar decomposition is
      $a=v|a|$, where $|a|$ is quarter--circular, $v$ is a Haar unitary
      and $\{v,|a|\}$ is $\ast$--free.
  \endproclaim

\noindent{\bf \S2\. Random matrices}

  The context for our study of random matrices will be the following.
  Let $\sigma$ be a probability measure on some measure space without atoms
    and let
    $L=\bigcap_{1\le p<\infty}L^{p}(\sigma)$
    (the algebra of complex valued random variables
    having all moments bounded) be
    endowed with the state $E$, given by $Ef=\int f\dif\sigma$.
  The $\ast$--algebra of $\nbyn$ random matrices is $M_{n}(L)=L\otimes\MnC$,
    with trace $\phi_n=E\otimes\tau_n$, where $\tau_n$ is the normalized
    trace on $\MnC$.
  For notation, let $\{e(i,j;n)\mid1\le i,j\le n\}$ be a system of
    matrix units in $\MnC$, and write for the $\nbyn$ random matrix,
    $A$, having $a_{ij}$ for $(i,j)$th entry
    $$ A = \sum_{1\le i,j\le n}a_{ij}e(i,j;n), \tag1 $$
    (i.e\. for convenience we omit the tensor product symbol).
  Thus we have
    $$ \phi_{n}(A)=n^{-1}\sum_{1\le k\le n}E(a_{kk}), $$
  We call $\Cpx\otimes\MnC\subseteq L\otimes\MnC$ the constant matrices,
    and we work with the block--diagonal constant matrices defined as follows.
  If $N$ is a positive integer dividing
    $n$, let $\Delta(N,n)\subseteq\Cpx\otimes\MnC$
    be $n/N$ copies of $\MNC$ embedded down the diagonal, i.e\.
    matrices of the form~(1) such that each $a_{ij}$
    is a constant and equals zero unless there are $0\le b\le(n/N)-1$ and
    $1\le\ih,\jh\le N$ such that $i=bN+\ih$ and $j=bN+\jh$.

  For an $\nbyn$ random matrix, $A=\sum_{i,j}a_{ij}e(i,j;n)$,
    denote by $E(A)$ the
    constant matrix whose $(i,j)$th entry is $E(a_{ij})$.
  Let
    $$ \align
      \brakphi{A}&=A-\phi(A)I, \\
      \brakE{A}&=A-E(A), \\
      \brakr{A}&=E(A)-\phi(A)I.
    \endalign $$

  \proclaim{Theorem 2.1}
    Let
      $$ Y(s,n)=\sum_{1\le i,j\le n}a(i,j;s,n)e(i,j;n) \in M_{n}(L) $$
      be random matrices for $s$ taking values in some set $S$
      such that
      \vskip1.5ex
      $$
      \left.
      \gathered
        a(i,j;s,n)=\overline{a(j,i;s,n)} \\ \vspace{ 0.5\jot}
      \aligned
        E( a(i,j;s,n))&=0 \\ \vspace{0.5\jot}
        E( |a(i,j;s,n)|^{2})&=\frac{1}{n}
      \endaligned
      \endgathered\;
      \right\}\qquad(1\le i\le j\le n,\,s\in S) \tag2
      $$
      \vskip1.5ex
      $$
        \sup_{1\le i\le j\le n}E(|a(i,j;s,n)|^{m})=O(n^{-m/2})
        \text{ as }n\rightarrow\infty\,\qquad(m\ge1,\,s\in S),
      \tag3
      $$
      that
      $$ \bigl\{\{a(i,j;s,n)\mid s\in S\}\mid1\le i\le j\le n\bigr\} \tag4 $$
      is an independent family of sets of random variables
      and that
      $$ E( a(i,j;s,n)a(j,i;s^{\prime},n))=0 \text{ whenever }
        s\neq s^{\prime}. \tag5 $$

    In addition, fixing $N$ and considering $n$ which are multiples of $N$,
      for each $t$ taking values in some set, $T$, let
      $$ D(t,n)=\sum \Sb 0\le b\le (n/N)-1\\ 1\le\ih,\jh\le N \endSb
        d(Nb+\ih,Nb+\jh;t,n)e(Nb+\ih,Nb+\jh;n) $$
      be elements of $\Delta(N,n)$,
      having a limit distribution as $n\rightarrow\infty$ and such that
      $$ \sup_{n,i,j}|d(i,j;t,n)|<\infty. $$
    Assume also that
      for each $t_{1},t_{2}\in T$,
      there exists $t_{3}\in T$ such that $D(t_{1},n)D(t_{2},n)=D(t_{3},n)$
      $\forall n$.

    Then the family of sets of random variables
      $$ \bigl(( \{Y(s,n)\})_{s\in S},\{ D(t,n)|t\in T\}\bigr) $$
      is asymptotically free as $n\rightarrow\infty$ passing through
      multiples of $N$, and moreover each
      $Y(s,n)$ has for limit distribution Wigner's semicircle law
      (of radius $2$).
  \endproclaim

  \demo{Proof}\footnote{Another proof of Gaussian cases is described
                        in the appendix.}
    We may assume without loss of generality
      that for some $t$, $D(t,n)=I_{n}$ $\forall n$,
      and that
      for each $t\in T$ and convergent sequences $\{\gamma_n\}_1^\infty$
      and $\{\lambda_n\}_1^\infty$ in $\Cpx$,
      there exists $t^{\prime}\in T$ such that
      $\gamma_nD(t,n)-\lambda_nI_{n}=D(t^{\prime},n)$
      $\forall n$.
    Hence to show asymptotic freeness, it suffices to show that if
      $s_{1},\ldots,s_{r}\in S$, $p_{1},\ldots,p_{r}\ge1$ and
      $t_{1},\ldots,t_{r}\in T$ are such that for each $1\le\iota\le r$,
      either $\phi_{n}(D(t_{\iota},n))=0\;\forall n$ or
      $D(t_{\iota},n)=I_{n}\;\forall n$, and in the latter case,
      $s_{\iota-1}\neq s_{\iota}$ (by convention $s_{0}=s_{r}$), then
      $$ \phi_{n}( D(t_{1},n)\brakphi{Y(s_{1},n)^{p_{1}}}\cdots
        D(t_{r},n)\brakphi{Y(s_{r},n)^{p_{r}}}) \tag6 $$
      approaches zero as $n\rightarrow\infty$.

    The proof that~(6) tends to zero is combinatorial
    and the heart of it is contained in Lemmas~2.2, 2.5 and~2.6.
    We shall see that a
      sharp enough counting argument is possible to show that~(6)
      goes to zero if we replace each
      $\brakphi{Y(s_{\iota},n)^{p_{\iota}}}$ by
      $\brakE{Y(s_{\iota},n)^{p_{\iota}}}$.
    Hence we need also show that the remainder
      $\brakr{Y(s_{\iota},n)^{p_{\iota}}}$
      goes to zero fast enough.
    Formally, we write
      $$ \brakphi{Y(s_{\iota},n)^{p_{\iota}}}
          =\brakE{Y(s_{\iota},n)^{p_{\iota}}}
          +\brakr{Y(s_{\iota},n)^{p_{\iota}}}, \tag7 $$
      substitute~(7) into~(6), expand
      and apply Lemmas~2.5 and~2.6.

    Lemma~2.3 shows that the $Y(s,n)$
      all have the same limit distribution, call it $\mu$.
    To see that $\mu$ is a semicircle law, one can prove it directly
      (see Remark~2.4);
     otherwise, one can use the central limit
      theorem as follows.
    Taking $S$ equal to the set, $\Naturals$, of natural numbers,
      for $M\ge1$, let
      $$ \psi(M,n)=M^{-1/2}\sum_{1\le s\le M}Y(s,n). $$
    Then $\psi(M,n)$ for $n\in\Naturals$ is a sequence of random matrices
      of the same type as the $Y(s,n)$'s
      and thus has
      limit distribution, $\mu$.
    But by the central limit theorem for free random variables (see~\cite{8}),
      the distributions of the $\psi(M,n)$'s
      tend to a semicircle law as $M\rightarrow\infty$.
    Hence $\mu$ is a semicircle law.
  \QED

  \proclaim{Lemma 2.2}
    Letting $Y(s,n)$ be as in Theorem~2.1 and
      fixing $s\in S$ and $m\ge1$, write
      $$ \brakr{Y(s,n)^{m}}=\sum_{1\le i,j\le n}f(i,j;n)e(i,j;n). $$
    Then
      $$ \align
        \sup_{1\le i<j\le n}|f(i,j;n)|&= O(n^{-3/2}) \\ \vspace{0.5\jot}
        \sup_{1\le k\le n}|f(k,k;n)|&=O(n^{-1/2})
      \endalign $$
      as $n\rightarrow\infty$.
  \endproclaim

  \demo{Proof}
    Let $1\le x,y\le n$.
    Then the $(x,y)$th element of $E(Y(s,n)^{m})$ is
      $$
        \bigl(E(Y(s,n)^{m})\bigr)_{xy}
        =\sum \Sb i_{2},\ldots,i_{m} \\ j_{k}=i_{k+1} \endSb
        E( a(x,j_{1};s,n)a(i_{2},j_{2};s,n)\cdots
        a(i_{m-1},j_{m-1};s,n) a(i_{m},y;s,n)), \tag8
      $$
      where the sum is taken over all $i_{2},\ldots,i_{m}\in\{1,\ldots,n\}$
      and $j_{k}$ is assigned the value $i_{k+1}$ for $1\le k\le m-1$.
    For ease of notation, we will frequently denote $x$ by $i_{1}$ and
      $y$ by $j_{m}$.

    In order to evaluate an expression of the form
      $$
        E( a(i_{1},j_{1};s,n)a(i_{2},j_{2};s,n)\cdots a(i_{m-1},j_{m-1};s,n)
        a(i_{m},j_{m};s,n))
      \tag9
      $$
      for a given choice of $i_1,\ldots,i_m$ and $j_1,\ldots,j_m$,
      one uses the independence condition,~(4),
      and
      gathers the $a(i_{k},j_{k};s,n)$'s together into mutually independent
      groups corresponding to the different (unordered) sets $\{i_{k},j_{k}\}$
      which appear.
    Then one multiplies together
      the expectations of each of the groups.
    Thus given a choice of
      $i_{1},\ldots,i_{m},j_{1},\ldots,j_{m}\in\{1,\ldots,n\}$,
      we define
      the resulting {\it gathering}, $\Qg$, to be the pair, $(\Rl,\Orient)$,
      where
      $\Rl$ is the equivalence relation on $\{1,\ldots,m\}$ indicating
      when $a(i_{k},j_{k};s,n)$ and $a(i_{ k^\prime},j_{ k^\prime};s,n)$
      belong to the same group, and $\Orient$ (the ``orientation map''
      of the gathering)
      indicates whether
      $a(i_{k},j_{k};s,n)$ and $a(i_{ k^\prime},j_{ k^\prime};s,n)$
      are on the same or opposite sides of the diagonal.
    Specifically, we let
      $$ k\Rtilde  k^\prime \text{ if and only if either }\;
        \text{(I)}\;\gathered  i_{k}=i_{ k^\prime} \\ \vspace{-2ex}
        \text{\scriptsize  and } \\ \vspace{-2ex}
        j_{k}=j_{ k^\prime} \endgathered
        \text{ or }\;\text{(II)}\;
        \gathered  i_{k}=j_{ k^\prime} \\ \vspace{-2ex}
        \text{\scriptsize  and } \\ \vspace{-2ex}
        j_{k}=i_{ k^\prime}, \endgathered \tag10 $$
      and define
      $$ \Orient:\bigl\{(k, k^\prime)\in\{1,\ldots,m\}^{2}\mid
        k\Rtilde k^\prime\bigr\}\rightarrow\{\pm1\} $$
      to be
      $$ \Orient(k, k^\prime)=\cases
          +1 & \text{ if (I)} \\
          -1 & \text{ if (II) but not (I).}
        \endcases \tag11
      $$
    We say that a gathering, $\Qg$, has property P if each
      equivalence class of $\Rl$ has at least two elements, and note that
      the quantity~(9) is nonzero only when its gathering has
      P, (because each $a(i,j;s,n)$ has zero expectation).

    Thus the sum in equation~(8) becomes
      $$
        \bigl(E(Y(s,n)^{m})\bigr)_{xy}
        =\sum_{ \Qg\text{ \tiny with \scriptsize P} }\quad
        \sum \Sb i_{2},\ldots,i_{m} \\
          j_k=i_{k+1} \\ \text{\tiny giving }\Qg \endSb
        E( a(x,j_{1};s,n)\cdots a(i_{m},y;s,n)),
      \tag12
      $$
      where the sums are over all possible gatherings, $\Qg$ with property P
      and all choices
      of $i_{2},\ldots,i_{m}\in\{1,\ldots,n\}$ which result in gathering $\Qg$.
    (We set $j_{k}=i_{k+1}$ for $1\le k\le m-1$.)

    Let us say that a gathering, $\Qg$
      has property P$_{3}$
      if it has property P and
      its equivalence relation, $\Rl$,
      also has an equivalence class of at least
      three elements.
    It has
      property P$_{4}$ if it
      has property P but not P$_{3}$ and also
      its orientation map, $\Orient$, takes on the value
      $+1$ for some pair $(k, k^\prime)$, where $k\neq k^\prime$.
    Those that are left (i.e\. those gatherings with property P but
      not P$_{3}$ or P$_{4}$) we say have property P$_{5}$.
    The sum over $\Qg$ with P may be split into sums over $\Qg$
      with P$_{3}$, P$_{4}$ and P$_{5}$.
    We shall show that the first two go to zero fast enough
       as $n\rightarrow\infty$, and the sum over $\Qg$ with P$_{5}$ approaches
       the value $\phi_{n}(Y(s,n)^{m})$ when $x=y$.

    Using independence and the bound on moments,~(3), we have
      $|E(a(i_1,j_1;s,n)\cdots a(i_m,j_m;s,n)|=O(n^{-m/2})$
      independent of the choice of $i$'s and $j$'s.
    Hence
      for a term corresponding to a fixed $\Qg$ in equation~(12),
      $$
        \sup_{1\le x\le y\le n}\left|\sum
          \Sb i_{2},\ldots,i_{m} \\
          j_k=i_{k+1} \\ \text{\tiny giving }\Qg \endSb
        E( a(x,j_{1};s,n)\cdots a(i_{m},y;s,n))\right|
        =O(n^{-m/2})\theta_{xy}(\Qg),
      \tag13
      $$
      where $\theta_{xy}(\Qg)$ is the number of choices of
      $i_1,\ldots,i_m$ which give $\Qg$.
    (Recall we set $j_{k}=i_{k+1}$ for $1\le k\le m-1$).

    A bound for $\theta_{xy}(\Qg)$
      may be found
      by doing computations using Feynmann graphs (similarly to their use
      in~\cite{11}).
    Consider the straight line graph (Figure~1).
    Associate the $k$th edge with $a(i_k,j_k;s,n)$ and the $k$th vertex
      with $i_k=j_{k-1}$ for $2\le k\le m$,
      with $x=i_1$ for $k=1$ and
      with $y=j_{m}$ for $k=m+1$.
    Then
      $$ \theta_{xy}(\Qg)\le n^{d(\Qg)-l_{xy}}, \tag14 $$
        where
      $d(\Qg)$ is the number of vertices in the quotient graph, $G$, obtained
      from the straight line graph
      by identifying edges $[k,k+1]$ and $[ k^\prime, k^\prime+1]$ if $k\Rtilde
k^\prime$
      (i.e\. if $a(i_{k},j_k;s,n)$ and $a(i_{ k^\prime},j_{k^\prime};s,n)$ are
      gathered by $\Qg$ into the same group), with orientation preserved
      if $\Orient(k, k^\prime)=+1$ and reversed if $\Orient(k, k^\prime)=-1$.
    We are thus identifying vertices which are forced by $\Qg$ to be labeled
      with the same value of $\{1,\ldots,n\}$.

    More explicitly,
     $d(\Qg)$ is the number of equivalence classes in $\{1,\ldots,m+1\}$
      modulo the equivalence relation $\vtilde$, which is generated by
      the relations
      $$
      \aligned
        &\left\{ \aligned
          \scriptstyle k&\vtilde k^\prime \\ \vspace{-1.5ex}
          \scriptstyle k+1&\vtilde k^\prime+1 \endaligned
        \bigm| k\Rtilde k^\prime
        \text{ and }\Orient(k, k^\prime)=+1\right\}\cup \\ \vspace{1\jot}
        \cup&\left\{ \aligned
        \scriptstyle k&\vtilde k^\prime+1 \\ \vspace{-1.5ex}
        \scriptstyle k+1&\vtilde k^\prime \endaligned
        \bigm| k\Rtilde k^\prime
        \text{ and }\Orient(k, k^\prime)=-1\right\}.
      \endaligned
      \tag15
      $$
    The fact that any assignment of values to $i_1,\ldots,i_m$
      which gives $\Qg$ must be constant on the equivalence
      classes of $\vtilde$ justifies equation~(14),
      when we also note that we may
      subtract $l_{xy}$ from $d(\Qg)$
      in~(14) to account for the fact that $i_{1}=x$
      and $j_{m}=y$ are already assigned.
    We may subtract $2$ if $1\not\vtilde m$ and $1$ if $1\vtilde m$, (which
      would imply that $x=y$).
    Thus let
      $$ l_{xy}=\cases 2&\text{if }x\neq y \\
        1&\text{if }x=y.\endcases $$

    Suppose that $\Qg$ has property P$_{3}$.
    Then $G$ has at most
      $\frac{m-1}{2}$ edges, hence
      (as $G$ is connected) at most $\frac{m+1}{2}$ vertices,
      so
      $$ \theta_{xy}(\Qg)\le n^{\frac{m+1}{2}-l_{xy}}. $$

    Suppose $\Qg$ has property P$_{4}$.
    Then a more intricate analysis of
      the identifications is called for in order to bound $d(\Qg)$.
    $\Rl$ groups the edges into pairs.
    Referring to~(15), observe that
      when $k\Rtilde k^\prime$, the identification of the edges
      $[k,k+1]$ and $[ k^\prime, k^\prime+1]$ results in each of the vertices
      being paired with a vertex not equal to itself
      except if the edges lie
      next to each other and are identified with reverse orientation
      (i.e\. $|k- k^\prime|=1$ and $\Orient(k, k^\prime)=-1$).
    For example, taking $k=3$ and $ k^\prime=4$, such an identification
      of neighboring edges with reverse orientation
      results in the graph shown in Figure~2.
    We call this a click, and we also use this word as a verb.
    A click removes two vertices from the line and leaves a tail
      of one vertex hanging below.
    Now, if possible, we click again, clicking any two edges of the line
      which are presently neighbors and are identified by $\Qg$
      with reverse orientation.
    After $c$ clicks the graph will consist of a line with $m-2c+1$ vertices
      together with tails hanging below which have
      a total of $c$ additional vertices.
    After we've clicked as much as possible, there will still remain
      some edges above on the line, because of the requirement that
      at least one pairing be orientation preserving.
    And then each of the vertices remaining on the line
      will become identified with
      some vertex not equal to itself
      when we do the further identifications.
    Hence $G$ will have at most
      $\frac{m-2c+1}2+c=\frac{m+1}{2}$ vertices and
      $$ \theta_{xy}(\Qg)\le n^{\frac{m+1}{2}-l_{xy}}. $$

    So
      equation~(12) becomes
      $$ ( E(Y(s,n)^{m}))_{xy}
        = S_5(x,y;m,n) + O(n^{1/2})n^{-l_{xy}}, \tag16 $$
      where
      $$
        S_5(x,y;m,n)=\sum_{\Qg\text{ \tiny with \scriptsize P}_{5}}\quad
        \sum
          \Sb i_{2},\ldots,i_{m} \\
          j_k=i_{k+1} \\ \text{\tiny giving }\Qg \endSb
        E( a(x,j_{1};s,n)\cdots a(i_{m},y;s,n))
        =n^{-m/2}\sum_{\Qg\text{ \tiny with \scriptsize P}_{5}}
        \theta_{xy}(\Qg)
      \tag17
      $$
      and where $O(n^{1/2})$ is uniform with respect to $x$ and $y$.
    (The last equality of~(17) results
      from equation~(2), independence and the fact that
      $\Qg$ with P$_{5}$ groups the $a(i,j;s,n)$'s into
      complex conjugate pairs.)
    The value of $S_5(x,y;m,n)$ thus depends only
      on whether $x=y$ or $x\neq y$ and not on the particular
      values of $x$ and $y$.
    Its easy to see that
      $\theta_{xy}(\Qg)=0$ when $x\neq y$ (and $\Qg$ has P$_{5}$).
    (Indeed, an induction argument shows, for example, that
      if we let $M$ be the number
      of elements $(i,j)$ in the list
      $(x,j_{1}),(i_{2},j_{2}),\ldots,(i_{m-1},j_{m-1}),(i_{m},y)$
      for which exactly one of $i$ and $j$ is equal to $x$, then $M$ is odd.)
    Hence we may conclude that
      $$ \sup_{1\le x<y\le n}\left|( E(Y(s,n)^{m}))_{xy}\right|
        =O(n^{-3/2}). $$

    Let $S(m,n)=S_5(x,x;m,n)$.
    To complete the proof of Lemma~2.2, it suffices to show
      that
      $$ \phi_{n}(Y(s,n)^{m})-S(m,n)=O(n^{-1/2}) \tag18 $$
      as $n\rightarrow\infty$.
    However~(18) follows from~(16) once we recall that
      $$ \phi(Y(s,n)^m)=n^{-1}\sum_{x=1}^nE(Y(s,n)^m)_{xx}. $$
  \QED

  Now let us show that each $Y(s,n)$ has a limit distribution.
  \proclaim{Lemma 2.3}
    For $Y(s,n)$ and in Theorem~2.1,
      for each $m\ge1$ there is an integer, $\alpha_{m}$,
      such that for every $s$,
      $$ \phi_{n}(Y(s,n)^{m})=\alpha_{m}+O(n^{-1/2}) $$
      as $n\rightarrow\infty$.
  \endproclaim

  \demo{Proof}
    We saw in Lemma~2.2 that
      $\phi_{n}(Y(s,n)^{m})=S(m,n)+O(n^{-1/2})$,
      so it suffices to show
      $$ S(m,n)=\alpha_{m}+O(n^{-1/2}), \tag19 $$
      where
      $$ S(m,n)=
        n^{-m/2}\sum_{\Qg\text{ \tiny with \scriptsize P}_{5}}
        \theta_{xx}(\Qg),
      \tag20
      $$
      and $\theta_{xx}(\Qg)$ is the cardinality of the set,
      $W_{\Qg}$, of choices of $i_{2},\ldots,i_{m}$ which give $\Qg$,
      (which, as we saw, is independent of $x\in\{1,\ldots,n\}$).
    For fixed $\Qg$ (with~P$_5$),
      let $G$ be the quotient graph obtained in the proof of Lemma~2.2.
    Either one or two vertices of $G$ are automatically assigned the
      value $x$, and there remain $d^{\prime}(\Qg)$ vertices of $G$ to which
      to assign values, (where $d^{\prime}(\Qg)=d(\Qg)-1$ or $d(\Qg)-2$).
    Then $W_{\Qg}$
      clearly contains those choices which assign distinct elements
      of $\{1,\ldots,n\}$ to these vertices of $G$, and is contained in
      the set of choices which assign arbitrary elements of $\{1,\ldots,n\}$
      to these vertices of $G$.
    Thus we see that
    $n(n-1)\cdots(n-d^{\prime}(\Qg))\le\theta_{xx}(\Qg)\le
n^{d^{\prime}(\Qg)}$.
    But the difference between these two bounds is $O(n^{d^{\prime}(\Qg)-1})$
      and $d^{\prime}(\Qg)\le m/2$ (because $G$ has $m/2$ edges).
    Hence~(19) holds if we set $\alpha_m$ equal to the number of
      gatherings, $\Qg$ with P$_5$, for which
      $d^\prime(\Qg)=m/2$.
  \QED

  \proclaim{Remark 2.4
    \footnote{Thanks to Alexandru Nica for showing me this.}}
    \rm
    One can now show that the distribution $\mu(X^m)=\alpha_m$ is
      a semicircle law of radius~$2$.
    Clearly $\alpha_m=0$ if $m$ is odd, so let $m=2h$ be even.
    Moreover, $\alpha_m$ is equal to the
      cardinality of the set, $\Cc$, of pairings of all the edges
      (with opposite orientation) of the straight line
      graph of length~$m$, such that successive clicking suffices to identify
      all the pairs.
    For each such pairing, $\Rl$, label each edge of the
      length~$m$ straight line graph with $+1$ or $-1$ as follows:
     for each $k\Rtilde k^\prime$ with $k< k^\prime$, label the $k$th edge with
      $-1$ and the $ k^\prime$th edge with $+1$.
    This labeling, which is actually a mapping from $\{1,\ldots,m\}$
      to $\{\pm1\}$, we call $\phi(\Rl)$.
    Let $\Ec$ be the set of all mappings
      $L:\{1,\ldots,m\}\rightarrow\{\pm1\}$ such that
      $\sum_{j=1}^mL(j)=0$ and $\sum_{j=1}^lL(j)\le0$ $\forall1\le l\le m$.

    $\phi$ maps $\Cc$ into $\Ec$, and it is actually a bijection.
    To see this, construct the inverse, $\psi$, of $\phi$,
     recursively as follows:
     if $L$ is a labeling, let $k$ be least such that $L(k)=+1$;
     then necessarily $L(k-1)=-1$;
     set $k-1\overset\psi(L)\to\sim k$, remove the edges $k-1$ and $k$
      and repeat the process, searching for the lowest pair
      of neighboring $-1$ and $+1$'s, until all the edges have been paired.

    Now, setting about determining the size of
      $\Ec$, it is easily seen that the set
      $\Ec$ is in bijection with the set of paths of length $m$ from
      $(0,0)$ to $(m/2,m/2)$ in the lattice $\{0,\ldots,m/2\}$ which
      never go below the diagonal.
    For example, if $m=6$, to the labeling $(-1,1,-1,-1,1,1)$
      we associate the path shown in Figure~3.
    The number of such paths can be computed without much difficulty
      (see for instance~\cite{16}, problem~81)
      and is found to be $\binom{2h}h-\binom{2h}{h-1}$.
    Hence it suffices to show that
      $$ \int_{-2}^2t^m\sqrt{4-t^2}\dif t
        =\cases
        0,&m\text{ odd} \\
        2\pi\left(\binom{2h}h-\binom{2h}{h-1}\right),&m=2h\ge2,\,
          m\text{ even},
        \endcases
      $$
      which is easily proved by induction.
  \endproclaim

  \proclaim{Lemma 2.5}
    Let $Y(s,n)$ and $D(t,n)$
      be as in Theorem~2.1.
    Let $t_{1},\ldots,t_{r}\in T$ and suppose that
      $D(t_{\iota},n)$ for $1\le\iota\le r$ is either equal to
      $I_{n}\;\forall n$ or has zero trace $\forall n$.
    Let $s_{1},\ldots,s_{r}\in S$ be such that if $D(t_{\iota},n)=I_{n}$
      then $s_{\iota-1}\neq s_{\iota}$, (where by convention $s_{0}=s_{r}$).
    Then
    $$
      \phi_{n}( D(t_{1},n)\brakE{Y(s_{1},n)^{p_{1}}}\cdots
      D(t_{r},n)\brakE{Y(s_{r},n)^{p_{r}}}) = O(n^{-1/2})
    \tag21
    $$
    as $n\rightarrow\infty$.
  \endproclaim

  \demo{Proof}
    Let
      $$ \align
        u_{\iota} &= p_{1}+\cdots+p_{\iota-1}+1\text{ for }1\le\iota\le r, \\
        m&=p_{1}+\cdots+p_{r}.
      \endalign
      $$
    Let us use the abbreviation
      $$
        C_{\iota}\text{ means }a(i_{u_{\iota}},j_{u_{\iota}};s_{\iota},n)
        a(i_{u_{\iota}+1},j_{u_{\iota}+1};s_{\iota},n)\cdots
        a(i_{u_{\iota}+p_{\iota}-1},j_{u_{\iota}+p_{\iota}-1};s_{\iota},n).
      \tag22
      $$
    Each $D(t,n)$ is the sum of its diagonal and off--diagonal parts, hence
      the left--hand side of equation~(21) is
      equal to a sum of $2^{r}$ terms of the form
      $$ \phi_{n}( \psi_{1}(t_{1},n)\brakE{Y(s_{1},n)^{p_{1}}}\cdots
      \psi_{r}(t_{r},n)\brakE{Y(s_{r},n)^{p_{r}}}),                   $$
      where each $\psi_{\iota}(t_{\iota},n)$ is either the diagonal or
      off--diagonal part of $D(t_{\iota},n)$, the same for all $n$.
    Thus
      $$ \aligned
        \phi_{n}(\psi_{1}(t_{1},n)&\brakE{Y(s_{1},n)^{p_{1}}}\cdots
          \psi_{r}(t_{r},n)\brakE{Y(s_{r},n)^{p_{r}}})         \\
\vspace{1\jot}
        &=n^{-1}
          \sum \Sb i_{1},\ldots,i_{m} \\ j_{k}=\ast \endSb
          \left(\prod_{1\le\iota\le r}
          d(j_{u_{\iota}-1},i_{u_{\iota}};t_{\iota},n)\right)
          E\bigl(( C_{1}-E(C_{1}))\cdots( C_{r}-E(C_{r}))\bigr)
          \\ \vspace{1\jot}
        &=n^{-1}
          \sum \Sb i_{1},\ldots,i_{m} \\ j_{k}=\ast \endSb
          \left(\prod_{1\le\iota\le r}
          d(j_{u_{\iota}-1},i_{u_{\iota}};t_{\iota},n)\right)
        \sum_{f\in\{0,1\}^{r}}(-1)^{|f|}
        E( E^{f_{1}}(C_{1})\cdots E^{f_{r}}(C_{r})),
      \endaligned \tag23 $$
      where we are summing
      over all $i_{1},\ldots,i_{m}$ in $\{1,\ldots,n\}$,
      and all allowable $j_{k}$, where $j_{k}$ is allowed to take on only
      the value $i_{k+1}$,
      unless $k=u_{\iota}-1$ and $\psi_{\iota}(t_{\iota},n)$ is off--diagonal,
      in which case we find $0\le b\le (n/N)-1$ and $1\le\ih\le N$ such that
      $i_{u_{\iota}}=Nb+\ih$ and allow $j_{k}$ to take on
      all of the values $Nb+\jh$ for $1\le\jh\le N$,
      $\jh\neq\ih$ (all subscripts of $i$ and $j$ are to be taken
      modulo $m$);
      we are summing over all $f=(f_{1},\ldots,f_{r})\in\{0,1\}^{r}$,
      and define $|f|=f_{1}+\cdots+f_{r}$,
      $$ E^{f_{\iota}}(C_{\iota})=\cases C_{\iota}&\text{if }
        f_{\iota}=0\\E(C_{\iota})&\text{if }f_{\iota}=1.\endcases $$

    Given a choice of $i_{1},\ldots,i_{m}$ and $j_{1},\ldots,j_{m}$,
      we define its gathering, $\Qg$, as in the proof of
      Lemma~2.2, (i.e\. $\Rl$ is defined by~(10) and $\Orient$ by~(11)).
    Note that once again, for a choice of $i_{1},\ldots,i_{m}$ and
      $j_{1},\ldots,j_{m}$
      to give a nonzero term
      in equation~(23), its gathering, $\Qg$, must
      have property P.
    Thus we may rearrange the sum to give
      $$ \aligned
        \phi_{n}&(\psi_{1}(t_{1},n)\brakE{Y(s_{1},n)^{p_{1}}}\cdots
          \psi_{r}(t_{r},n)\brakE{Y(s_{r},n)^{p_{r}}}) \\ \vspace{0.5\jot}
        &=\sum_{\Qg\text{ \tiny with \scriptsize P}}n^{-1}\sum
            \Sb i_{2},\ldots,i_{m} \\ j_k=\ast \\
            \text{\tiny giving }\Qg \endSb
          \left(\prod_{1\le\iota\le r}
          d(j_{u_{\iota}-1},i_{u_{\iota}};t_{\iota},n)\right)
          \sum_{f\in\{0,1\}^{r}}(-1)^{|f|}
          E( E^{f_{1}}(C_{1})\cdots E^{f_{r}}(C_{r})),
      \endaligned
      \tag24
      $$
      where the sums are over all
      gatherings, $\Qg$, with P
      and all choices of $i_{1},\ldots,i_{m}$
      and all allowable $j_{1},\ldots,j_{m}$
      in $\{1,\ldots,n\}$ giving gathering $\Qg$.
    Let $\theta(\Qg)$ be the number of such choices.
    Using the generalized H\"older inequality together with
      the boundedness conditions on the moments of the $a(i,j)$'s and
      on the absolute values of the $d$`s,
      we have that for a particular $\Qg$, the corresponding term in
      equation~(24) becomes
        $$
          n^{-1}\sum
            \Sb i_{2},\ldots,i_{m} \\ j_k=\ast \\
            \text{\tiny giving }\Qg \endSb
            \left(\prod_{1\le\iota\le r}
            d(j_{u_{\iota}-1},i_{u_{\iota}};t_{\iota},n)\right)
            \sum_{f\in\{0,1\}^{r}}(-1)^{|f|}
            E( E^{f_{1}}(C_{1})\cdots E^{f_{r}}(C_{r}))
            =O(n^{-\frac{m}{2}-1})\theta(\Qg).
        \tag25
        $$
    We will show that for each $\Qg$ with P, the
      left hand side of~(25) is $O(n^{-1/2})$ as
      $n\rightarrow\infty$.

    For each $\iota$ with $\psi_{\iota}(t_{\iota},n)$ off--diagonal
      (and $i_{u_{\iota}}=Nb+\ih_{\iota}$), there are $N-1$ allowable values of
      $j_{u_{\iota}-1}=Nb+\jh_{\iota}$ all nonequal to
      $i_{u_{\iota}}$ and corresponding to choices of
$\ih_{\iota}-\jh_{\iota}$.
    Hence
      $$ \theta(\Qg)=\sum
          \Sb h_{\iota}\in\{1,\ldots,N-1\} \\
          \text{\scriptsize for }
          \psi_{\iota}\text{\scriptsize off--diag.}\endSb
        \theta(\Qg,h), \tag26 $$
      where each $h_\iota$ specifies the value of $\ih_\iota-\jh_\iota$
      modulo $N$, and we denote by $h$ the aggregate of the $h_\iota$'s, for
      $\iota$ running over those values for which $\psi_\iota$
      is off--diagonal;
     $\theta(\Qg,h)$ is the number of choices of $i_1,\ldots,i_m$
      such that when $j_1,\ldots,j_m$ are then taken
      according to $h$, $\Qg$ is the resulting gathering;
     in~(26), we are summing over all possible values of
      $h=(h_\iota)_{\{\iota\mid\psi_\iota\text{ off--diag.}\}}$.
    Fixing gathering $\Qg=(\Rl,\Orient)$, let $G$ be the quotient graph
      obtained from the $m$--gon graph (see Figure~4) by, as in the
      proof of Lemma~2.2, identifying edges according to $\Rl$ with
      orientations specified by $\Orient$,
      and let $d(\Qg)$ be the number of vertices of $G$.
    So $d(\Qg)$ is the number of equivalence classes of $\{1,\ldots,m\}$
      under $\vtilde$, defined at~(15).

    We claim that
      $$ \theta(\Qg,h)\le n^{d(\Qg)} \tag27 $$
      for every $h$.
    To see this, note that for $k_1\vtilde k_2$, $i_{k_1}$ and $i_{k_2}$
      differ by a fixed integer (determined by $h$).
    Thus there are $d(\Qg)$ degrees of freedom when choosing values of
      $i_1,\ldots,i_m$ giving gathering $\Qg$,
      which proves~(27).

    As an aside, let us remark at this point that for a given choice
      of $\Qg$ and $h$, the above mentioned fixed differences between
      $i_{k_1}$ and $i_{k_2}$ for $k_1\vtilde k_2$ my be incompatible
      within an equivalence class of $\vtilde$, in which case
      $\theta(\Qg,h)=0$.
    This phenomena will be important later.

    The case of $\Qg$ with property~P
      may be split into the three cases, $\Qg$ with
      P$_{3}$, P$_{4}$ and P$_{5}$.
    If $\Qg$ has P$_{3}$, then the quotient graph, $G$ has at most
      $\frac{m-1}{2}$ edges, so $\frac{m+1}{2}$ vertices.
    Hence $d(\Qg)\le\frac{m+1}{2}$ and the quantity,~(25),
      is $O(n^{-1/2})$ as $n\rightarrow\infty$.

    If $\Qg$ has P$_{4}$, then one may click $c$ times to
      leave a ring with $m-2c$ vertices and tails with
      $c$ vertices.
    (See the proof of Lemma~2.2 for the definition of a click,
      or just look at Figure~2.)
    If one has clicked as much as possible and is still left with
      an inner ring, then each vertex in the inner ring becomes identified
      with some vertex not equal to itself
      after the remaining edge identifications are
      made, so $G$ has at most $\frac{m}{2}$ vertices.
    If one is able to click so much that the ring disappears, then since
      by hypothesis
      one pair of vertices is identified with orientation preserved
      (and thus can never take part in a click), the only possibility is,
      before making the last click, to have a square with tails.
    Now clicking two of the sides of the square
      and identifying the
      other two with orientation preserved has the effect of identifying
      three vertices of the square together.
    Thus $G$ has at most $\frac{m}{2}$ vertices.
    Hence the quantity,~(25),
      is $O(n^{-1})$ as $n\rightarrow\infty$.

    If $\Qg$ has P$_{5}$, then $\Qg$ pairs each $a(i,j;s_{\iota},n)$
      with some $a(j,i;s_{\iota^{\prime}},n)$.
    We may suppose that always $s_{\iota^{\prime}}=s_{\iota}$,
      for otherwise the left hand side of~(25)
      is zero by condition~(4).
    Moreover $E(C_{\iota})$ is zero unless $\Qg$ pairs $a(i,j;s_{\iota},n)$'s
      only with other $a(i,j;s_{\iota},n)$'s from the same $C_{\iota}$.
    Such a $\Qg$ we say preserves $C_{\iota}$.
    Let us say that a given $f\in\{0,1\}^{r}$ is compatible with
      $\Qg$ if $f_{\iota}=1$ implies that $\Qg$ preserves $C_{\iota}$.
    Then
      $$ E( E^{f_{1}}(C_{1})\cdots E^{f_{r}}(C_{r}) )=
        \cases
          n^{-m/2}&\text{if $f$ is compatible with $\Qg$} \\
          0&\text{otherwise.}\endcases $$
    Now it is easy to show that
      $$
        \sum_{f\in\{0,1\}^{r}}(-1)^{|f|}
        E( E^{f_{1}}(C_{1})\cdots E^{f_{r}}(C_{r}) )
        = \cases
          n^{-m/2}
          &\text{if $\Qg$ preserves no $C_{\iota}$}  \\
          0 & \text{otherwise.}
          \endcases
      \tag28
      $$

    Thus we are left to show that the
      the left hand side of~(25)
      is $O(n^{-1/2})$ as
      $n\rightarrow\infty$ for $\Qg$ which has P$_{5}$, pairs
      each random variable
      $a(i,j;s_{\iota},n)$ with its conjugate,
      and preserves no $C_{\iota}$.
    Because $\Qg$ doesn't preserve any $C_{\iota}$, there
      must be at least one edge in each $C_{\iota}$ which is paired with
      an edge in some $C_{\iota^{\prime}}$, for $s_{\iota}=s_{\iota^{\prime}}$
      and $\iota\neq\iota^{\prime}$.
    Suppose first that each $D(t_{\iota},n)=I_{n}$ and hence that
      $s_{\iota-1}\neq s_{\iota}$ for every $\iota$.
    Then between $C_{\iota}$ and $C_{\iota^{\prime}}$ on the graph,
      there are always
      edges belonging to some $C_{\kappa}$,
      where $s_{\iota}\neq s_{\kappa}$.
    Hence by clicking, we'll never be able to completely remove all the
      edges of any $C_{\iota}$ from the ring.
    Thus, as we've seen, clicking as much as possible ($c$ times) will leave
      $m-2c$ vertices in the ring and $c$ in the tails, and each of
      the vertices remaining in the ring will become identified with at
      least one other vertex not itself.
    $G$ will thus have at most $m/2$ vertices, and the
      the left hand side of~(25)
      is $O(n^{-1})$ as $n\rightarrow\infty$.

    Now suppose that some $D(t_{\kappa},n)$ has trace zero.
    We saw in the last paragraph that
      the quantity~(25) is $O(n^{-1})$
      unless we suppose that by merely clicking, one may eventually identify
      an edge of some $C_{\iota}$ with an edge of $C_{\iota-1}$.
    But because $\Qg$ preserves no $C_{\kappa}$, for the first such $\iota$,
      call it $\iotat$,
      this would imply that the vertex numbered $u_{\iotat}$
      gets clicked
      out of the ring into the tails without being identified with any
      of the vertices,
      $\{u_{\iota}\mid \iota\neq\iotat\}$.
    Clearly $s_{\iotat}=s_{\iotat-1}$ so $D(t_{\iotat},n)$ must have
      trace zero, and we claim that in addition,
      $\psi_{\iotat}(t_{\iotat},n)$
      is the diagonal part of $D(t_{\iotat},n)$.
    (This is where the aforementioned phenomena of incompatibility
      of $\Qg$ and $h$ comes into the limelight.)
    Clicking may have occurred within $C_{\iotat-1}$ and
      $C_{\iotat}$, and then a click occurs at the boundary between them.
    Thus, for this boundary--click, an edge of $C_{\iotat-1}$, call it
      $a(i_{k_{1}},j_{k_{1}};s_{\iotat-1},n)$, is clicked with an
      edge of $C_{\iotat}$, call it $a(i_{k_{2}},j_{k_{2}};s_{\iotat},n)$;
      moreover, tracing the history of clicks, we see that
      $j_{k_{1}}=j_{u_{\iotat}-1}$ and
      $i_{k_{2}}=i_{u_{\iotat}}$.
    Then, as $\Orient(k_{1},k_{2})=-1$, we have
      $j_{u_{\iotat}-1}=j_{k_{1}}=i_{k_{2}}=i_{u_{\iotat}}$,
      which implies that $\psi_{\iotat}(t_{\iotat},n)$ is diagonal.

    If $h$ is specified and
      we know the quotient graph $G$, what do we know about the set
      of choices of $i_{1},\ldots,i_{m}$ which give $\Qg$?
    It clearly contains those choices which assign to the vertices of $G$
      elements of $\{1,\ldots,n\}$ which are pair--wise of distance
      at least $w$ apart, and is contained in
      the set of choices which assign arbitrary elements of $\{1,\ldots,n\}$
      to the vertices of $G$.
    Thus we see that $n(n-2w)\cdots(n-d(\Qg)(2w))\le\theta(\Qg)\le n^{d(\Qg)}$.
    But the difference between these two bounds is $O(n^{d(\Qg)-1})$, and
      as $d(\Qg)\le1+m/2$ (because $G$ has $m/2$ edges), it does no
      harm in equation~(25) to sum
      over those choices of $i_{1},\ldots,i_{r}$ which assign arbitrary
      elements of $\{1,\ldots,n\}$ to the vertices of $G$.
    However we see from~(28) that, in equation~(25),
      since $\Qg$ has P$_{5}$
      and preserves no $C_{\iota}$, each sum over $f\in\{0,1\}^{r}$ has value
      $n^{-m/2}$.
    Moreover as we saw in the preceding paragraph,
      the vertex $u_{\iotat}$ is identified with no other $u_{\iota}$
      vertex in $G$, so the quantity
      $\sum_{1\le k\le n}d(k,k;t_{\iotat},n)$ (which is zero because
      it equals the trace of $D(t_{\iotat},n)$)
      factors out of equation~(25).
    Thus even if one can click away the ring, the left hand side of~(25)
      equals zero plus $O(n^{-1})$.

    Hence we have proved that each of the terms of
      equation~(24) is $O(n^{-1/2})$.
  \QED

  \proclaim{Lemma 2.6}
    Let $Y(s,n)$ and $D(t,n)$ be as in Theorem~2.1.
    Let $s_{1},\ldots,s_{r}\in S$, $t_{1},\ldots,t_{r}\in T$.
    Then
    $$ \phi_{n}( D(t_1,n)\brakstar{Y(s_{1},n)^{p_{1}}}\cdots
      D(t_r,n)\brakstar{Y(s_{r},n)^{p_{r}}}) = O(n^{-1/2}) \tag29 $$
    as $n\rightarrow\infty$,
    where each $\brakstar{Y(s_{\iota},n)^{p_{\iota}}}$ means either
    $\brakE{Y(s_{\iota},n)^{p_{\iota}}}$
    or $\brakr{Y(s_{\iota},n)^{p_{\iota}}}$,
    the same for all $n$, and it is $\brakr{Y(s_{\iota},n)^{p_{\iota}}}$
    for at least one $\iota$.
  \endproclaim
  \demo{Proof}
    Let $l$ be the number of $\brakstar{Y(s_{\iota},n)^{p_{\iota}}}$ in
      equation~(29) which are $\brakr{Y(s_{\iota},n)^{p_{\iota}}}$,
      (thus $l\ge1$), and
      let $m$ be the sum of those $p_{\iota}$ for which
       $\brakstar{Y(s_{\iota},n)^{p_{\iota}}}$ is
       $\brakE{Y(s_{\iota},n)^{p_{\iota}}}$.
    Write
      $$ \brakr{Y(s_{\iota},n)^{p_{\iota}}}
        =\sum_{1\le i,j\le n}f(i,j;s_{\iota},p_{\iota},n)e(i,j;n), $$
    Then
      $$ \phi_{n}( D(1,n)\brakstar{Y(s_{1},n)^{p_{1}}}\cdots
        D(r,n)\brakstar{Y(s_{r},n)^{p_{r}}})
        =n^{-1}\sum \Sb i_{1},\ldots,i_{l+m} \\ j_{k}=\ast \endSb
        cE(\Pc),
      \tag30
      $$
      where $c$ (depending on the $i$'s and $j$'s) is a product of
      $d(j_{u_{\iota}-1},i_{u_{\iota}};t_{\iota},n)$'s, one for each $\iota$,
      $\Pc$ (also depending on the $i$'s and $j$'s)
      is a product of $l$ terms of the form
      $f(i_{k},j_{k};s_{\iota},p_{\iota},n)$
      and $r-l$ terms of the form
      $C_{\iota}-E(C_{\iota})$
      and the sum is over all $i_{1},\ldots,i_{m}$ and allowable
      $j_{1},\ldots,j_{m}$ in $\{1,\ldots,n\}$,
      where ``allowable'' is defined as in the proof of Lemma~2.5,
      modified in the obvious way for this situation.
    Let $L$ be the set of $k\in\{1,\ldots,l+m\}$ such that
      $f(i_{k},j_{k};p_{\iota},n)$
      appears in $\Pc$, and let $M$ be
      the compliment of $L$.
    Thus $L$ and $M$ have cardinalities $l$ and $m$, respectively.
    To find the precise expressions for $C_{\iota}$ as
      at~(22), we would
      need to redefine $u_{\iota}$, but the relevant fact here
      is that $C_{\iota}$ is a product of $p_{\iota}$ random variables of
      the form $a(i,j;s_{\iota},n)$;
     moreover, the redefinition of the $u_{\iota}$'s is necessary to
      correctly describe the rules of ``allowable'' $j_{k}$'s and to
      define $c$, but the relevant facts are that for
      $r$ values of $k$ there are at most $N$
      values $j_{k}-i_{k+1}$ can take,
      that $j_k=i_{k+1}$ for all other values of $k\in M$,
      and that $c$ is bounded.
    As in the proofs of the above lemmas, for each choice of
      $i_{1},\ldots,i_{m}$ and $j_{1},\ldots,j_{m}$,
      let $\Qg=(\Rl,\Orient)$ be the
      associated gathering on $M$, i.e\. $\Rl$ is an equivalence
      relation on $M$ indicating how the $a(i,j;n)$'s
      group together
      (defined precisely by~(10))
      and $\Orient$ indicates the orientations
      of the identifications, and is defined precisely by~(11).
    As before, for $E(\Pc)$ to be nonzero, $\Qg$ must have property P.
    Moreover, let $L_{d}$ having cardinality
      $l_{d}$ (``$d$'' for diagonal) be the set of
      $k\in L$ such that $i_{k}=j_{k}$ and let $L_{o}$
      having cardinality $l_{o}$ (``$o$'' for off--diagonal)
      be the set of
      $k\in L$ such that $i_{k}\neq j_{k}$.
    Then
      Lemma~2.2 and the boundedness of moments~(3), together with the
      generalized H\"older inequality imply that
      $$ E(\Pc) = O( n^{-(m+l_{d}+3l_{o})/2}) $$
      as $n\rightarrow\infty$, where $O(n^{-(\cdots)/2})$ is
      independent of $i_1,\ldots,i_{l+m}$ and $j_1,\ldots,j_{l+m}$.
    Thus rearranging the sum in equation~(30),
      $$ \aligned
        \phi_{n}( D(1,n)\brakstar{Y(s_{1},n)^{p_{1}}}\cdots
          D(r,n)\brakstar{Y(s_{r},n)^{p_{r}}})
          \;&=\;n^{-1}\sum_{\Qg\text{ \tiny with \scriptsize P}}\quad
            \sum_{L_{d}\subseteq L}\quad
            \sum \Sb i_{1},\ldots,i_{l+m}\\ j_{k}=\ast \endSb
            cE(\Pc)              \\ \vspace{1\jot}
          &= n^{-1}\sum_{\Qg\text{ \tiny with \scriptsize P}}\quad
            \sum_{L_{d}\subseteq L} O( n^{-(m+l_{d}+3l_{o})/2})
            \theta(\Qg,L_{d}),
        \endaligned
      \tag31
      $$
      where the sums are over all gatherings $\Qg$, on $M$ having property
      P, all $L_{d}$, subsets of $L$, and all $i_{1},\ldots,i_{l+m}$
      and all allowable $j_{1},\ldots,j_{l+m}$ in
      $\{1,\ldots,n\}$ which give gathering $\Qg$ and diagonal set $L_{d}$;
     and where $\theta(\Qg,L_{d})$ is the number of such choices of
      $i_{1},\ldots,i_{l+m}$ and $j_{1},\ldots,j_{l+m}$.

    Now
      $$ \theta(\Qg,L_{d})\le N^rn^{d(\Qg,L_{d})}, \tag32 $$
      (we take $N^r$ to account for different choices of $j_{1},\ldots,j_{m+l}$
      possible for each $i_{1},\ldots,i_{m+l}$),
      where $d(\Qg,L_{d})$ is the number of vertices in the
      quotient graph $G$ of the $(m+l)$--gon graph,
      where $G$ is obtained in the following
      way:
      \roster
        \item"I)" if $k\in L_{o}$, do nothing to the $k$th edge,
        \item"II)" if $k\in L_{d}$, remove the $k$th edge and identify its
          vertices,
        \item"III)" if $k, k^\prime\in M$ and $k\Rtilde k^\prime$, then
identify the
          $k$th and $ k^\prime$th edges with orientation determined by
          $\Orient(k, k^\prime)$.
      \endroster
    Thus $G$ has at most $l_o+(m/2)$ edges, so
      $d(\Qg,L_{d})\le l_{o}+\frac{m}{2}+1$.
    Hence substituting into equations~(32)
      and~(31), and noting that
      either $l_{d}\ge1$ or $l_{o}\ge1$,
      one has proved the lemma.
  \QED

  Thus the proof of Theorem~2.1
    is completed.
  A special case is
  \proclaim{Corollary 2.7}
    If $N$ divides $n$ let $\Psi_{Nn}$ embed $M_{N}(\Cpx)$ in $M_{n}(L)$
      by
      $$ \Psi_{Nn}( e(\ih,\jh;N))
        =\sum_{b=0}^{(n/N)-1}e(Nb+\ih,Nb+\jh;n) $$
      for $1\le\ih,\jh\le N$.
    Let $Y(s,n)$ be random matrices as in
      Theorem~2.1.
    Then the family of subsets of random variables,
      $$ \bigl(( Y(s,n))_{s\in S},
        \Psi_{Nn}(M_{N}(\Cpx))\bigr) $$
      is asymptotically free as $n\rightarrow\infty$ passing through
      multiples of $N$.
  \endproclaim

  As can be seen, in the proofs of Lemmas~2.5
    and~2.6, the net effect of $N$ is a constant multiple
    of $N^{r}$ in the bounds on the asymptotic behavior of the quantities
    under consideration.
  If $N$ were allowed to grow as $n\rightarrow\infty$, the only criteria
    for the theorem to remain true
    would be that $N^{r}O(n^{-1/2})$ go to zero as $n\rightarrow\infty$
    (for every $r$).
  Hence we may extend Theorem~2.1 to the
    following:
  \proclaim{Theorem 2.8}
    Let $Y(s,n)$ be as in Theorem~2.1.
    Let $n_{k}$ be a sequence of positive integers tending to $\infty$
      and let $N_{k}$ dividing $n_{k}$ for each $k$ be such that
      $N_{k}n_{k}^{-\epsilon}$ goes to zero as $k\rightarrow\infty$
      for all $\epsilon>0$.
    For each $t$ taking values in some set, $T$, let
      $$ D(t,n_{k})=\sum \Sb 0\le b\le (n_{k}/N_{k})-1 \\
          1\le\ih,\jh\le N_{k} \endSb
        d(N_{k}b+\ih,N_{k}b+\jh;t,n_{k})e(N_{k}b+\ih,N_{k}b+\jh;n_{k}) $$
      be elements of $\Delta(N_{k},n_{k})$
      having a limit distribution as $k\rightarrow\infty$ and such that
      $$ \sup_{k,i,j}|d(i,j;t,n_{k})|<\infty. \tag33 $$
    Assume also that
      for each $t_{1},t_{2}\in T$,
      there exists $t_{3}\in T$ such that
      $D(t_{1},n_{k})D(t_{2},n_{k})=D(t_{3},n_{k})$
      $\forall k$.
    Then the family of sets of random variables
      $$ \bigl(( \{Y(s,n_{k})\})_{s\in S},
        \{ D(t,n_{k})|t\in T\}\bigr) $$
      is asymptotically free as $k\rightarrow\infty$
      and moreover each
      $Y(s,n_{k})$ has limit distribution equal to a semicircle law.
  \endproclaim
  Also in the above, $N_{k}$ could be allowed to grow more quickly if also
    the bounds,~(33)
    on $d(i,j;t,n_{k})$ went to zero quickly enough as
    $k\rightarrow\infty$.

\noindent{\bf\S3\. Free products of $L(\freeF_K)$ with a matrix algebra}

  In this section we use random matrices to
    investigate  $L(\freeF_K)\ast M_N(\Cpx)$.
  Henceforth, as pertains to elements of $\ast$--noncommutative probability
    spaces, ``semicircular'' will mean semicircular of radius $2$,
    and similarly ``quarter--circular'' and ``circular'' will mean
    of radius $2$.
  We shall use the following result, which Voiculescu proved
    by
    approximating a semicircular element with random matrices and breaking
    the $\nbyn$ approximants into $N^{2}$ $(n/N\times n/N)$ blocks.

  \proclaim{Proposition 3.1}{\rm (\cite{10}).}
    Let $(A,\phi)$ be a $\ast$--noncommutative probability space and let
      $D\subseteq A$ be a commutative $\ast$--subalgebra.
    Let $A$ contain
      $$ \align
        \nu&=\{F(p;s)\mid1\le p\le N,s\in S\}
          \text{ --- a free family of semicircular elements,} \\
        \nu^{\prime}&=\{G(p,q;s)\mid1\le p<q\le N,s\in S\}
          \text{ --- a $\ast$--free family of circular elements,}
      \endalign
      $$
      such that $\{\nu,\nu^\prime,D\}$ is $\ast$--free.
    In $A\otimes M_{N}(\Cpx)$ let
      $$ \align
        X(s)&=N^{-1/2}\left(\sum_{1\le p\le N}F(p;s)\otimes e_{pp}
        + \sum_{1\le p<q\le N}( G(p,q;s)\otimes e_{pq}+
        G(p,q;s)^\ast\otimes e_{qp})\right) \\ \vspace{1\jot}
        d&=\sum_{1\le p\le N}d(p)\otimes e_{pp},
      \endalign
      $$
      where $\{d(p)\mid1\le p\le N\}\subseteq D$.
    Then each $X(s)$ is semicircular and $\{d,( X(s))_{s\in S}\}$
      is $\ast$--free in $( A\otimes M_{N}(\Cpx),\phi\otimes\tau_{N})$.
  \endproclaim

  The next result uses the random matrices with block diagonals
    of \S2.

  \proclaim{Proposition 3.2}
    Let $S$ be any nonempty set, $N\ge2$.
    In a $\ast$--noncommutative probability space $(A,\phi)$
      with $\phi$ a trace,
      let $\nu_{1}=\{X(s)\mid s\in S\}$ be a free family of
      semicircular elements and
      $\nu_{2}=\{e_{ij}\mid1\le i,j\le N\}$ a system of matrix units in $A$
      such that $\{\nu_{1},\nu_{2}\}$ is free.
    Let
      $$ \matrix\format\r&\l\qquad&\l\\
        F(k,s)&=N^{1/2}e_{1k}X(s)e_{k1}&(1\le k\le N) \\ \vspace{1\jot}
        G(i,j,s)&=N^{1/2}e_{1i}X(s)e_{j1}&(1\le i<j\le N)
      \endmatrix $$
      for $s\in S$, $1\le i,j,k\le N$, $i<j$.
    Then in the $\ast$--noncommutative probability space,
      $( e_{11}Ae_{11},N\phi\restrict_{e_{11}Ae_{11}})$,
      $$ \omega_{1}=\{F(k,s)\mid1\le k\le N,\:s\in S\} $$
      is a free family of
      semicircular elements and
      $$ \omega_{2}=\{G(i,j,s)\mid1\le i<j\le N,\:s\in S\} $$
      is a $\ast$--free family
      of circular elements such that $\{\omega_{1},\omega_{2}\}$ is
      $\ast$--free.
  \endproclaim

  \demo{Proof}
    Model the semicircular family $\nu_1$ and the matrix units
      as in Corollary~2.7, with the additional
      stipulation that $\{a(i,j;s,n)\mid1\le i\le j\le n,\,s\in S\}$
      be an independent set of random variables in $L$.
    Consider the families, $\Omega_1(n)$ and $\Omega_2(n)$
      of $(n/N\times n/N)$  matrices which approximate
      the families $\omega_1$ and $\omega_2$.
    Take $\sqrt2$ times the real and imaginary parts of
      the matrices in $\Omega_2(n)$ and apply Theorem~2.1 to conclude
      that they, together
      with $\Omega_1(n)$, are asymptotically free, each with
      semicircular limit distribution.
  \QED

  In a similar way we can prove

  \proclaim{Proposition 3.3}
    Let $S$ be any nonempty set, $N=2M\ge4$.
    In a $\ast$--noncommutative probability space $(A,\phi)$
      with $\phi$ a trace,
      let $\nu_{1}=\{X(s)\mid s\in S\}$ be a free family of
      semicircular elements and
      $\nu_{2}=\{e_{ij}\mid1\le i,j\le N\}$ a system of matrix units in $A$
      such that $\{\nu_{1},\nu_{2}\}$ is free.
    Consider the system of $2\times2$ matrix units
      $E_{11}=\sum_{k=1}^Me_{kk}$, $E_{12}=\sum_{k=1}^Me_{k,M+k}$,
      $E_{21}=E_{12}^*$ and $E_{22}=1-E_{11}$.
    Then in the $\ast$--noncommutative probability space,
      $( E_{11}AE_{11},2\phi\restrict_{E_{11}AE_{11}})$,
      the family
      $$ \bigl(\{E_{11}XE_{11}\},\{E_{12}XE_{21}\},\{E_{11}XE_{21}\},
        \{e_{pq}\mid1\le p,q\le M\}\bigr) $$
      is $*$--free.
  \endproclaim

  \proclaim{Theorem 3.4}
    $L(\freeF_{K})\ast M_{N}(\Cpx)
      \cong L(\freeF_{N^{2}K})\otimes M_{N}(\Cpx)$.
  \endproclaim

  \demo{Proof}
    By straightforward algebra, (e.g\. \cite{10}, Lemma~3.1), one can
      apply Proposition~3.2 to show that
      $$ (L(\freeF_K\ast M_N(\Cpx))_{1/N}\cong L(\freeF_{N^2K}). $$
  \QED

  Voiculescu proved~(\cite{10} Theorem~3.3)
    that
    $$
      L(\freeF_{N^{2}K+1})\otimes M_{N}(\Cpx)
      \cong L(\freeF_{K+1}),
    \tag35
    $$
    so we might hope to get something
    similar for $L(\freeF_{N^{2}K})\otimes M_{N}(\Cpx)$
    which ``approaches''~(35)
    as $N$ gets large.
  We indeed get such a result using the following
    proposition, which is little more than the proof
    of~\cite{10} Theorem~3.3
    run backwards.
  The idea is, by conjugating with unitaries, to show that Proposition~3.1
    holds also when we add a matrix like one of the $X(s)$'s,
    but with the circular $G(1,q)$'s and $G(q,1)$'s for
    $2\le q\le N$ replaced by quarter--circular elements.

  \proclaim{Proposition 3.5}
    Let $N$ be a positive integer, $S$ a set (finite or
      infinite, possibly empty).
    Let $(A,\phi)$ be a W$^*$--probability space with $\phi$ a trace,
      such that $A$ contains
      $$ \align
        \nu_{1}&=\{a(p)\mid1\le p\le N\}
          \text{ --- a $\ast$--free family of normal elements,} \\
        \nu_{2}&=\{f(p)\mid1\le p\le N\}
          \text{ --- a free family of semicircular elements,} \\
        \nu_{3}&=\{g(p,q)\mid2\le p<q\le N\}
          \text{ --- a $\ast$--free family of circular elements,} \\
        \nu_{4}&=\{b(q)\mid2\le q\le N\}
          \text{ --- a free family of quarter--circular elements,} \\
        \nu_{5}&=\{f(p;s)\mid1\le p\le N,s\in S\}
          \text{ --- a free family of semicircular elements,} \\
        \nu_{6}&=\{g(p,q;s)\mid1\le p<q\le N,s\in S\}
          \text{ --- a $\ast$--free family of circular elements,}
      \endalign
      $$
      such that $\{\nu_{1},\nu_{2},\nu_{3},\nu_{4},
      \nu_{5},\nu_{6}\}$ is $\ast$--free.
    Let $M_{N}(\Cpx)$ have system of matrix units $\{e_{pq}\mid1\le p,q\le
N\}$.
    In $(A\otimes M_{N}(\Cpx),\phi\otimes\tau_N)$, let
      $$ \align
        H&=\sum_{1\le p\le N}f(p)\otimes e_{pp}
          + \sum_{2\le q\le N}b(q)\otimes(e_{1q}+e_{q1})
          + \sum_{2\le p<q\le N}( g(p,q)\otimes e_{pq}+
          g(p,q)^\ast\otimes e_{qp}) \\ \vspace{1\jot}
        C&=\sum_{1\le p\le N}a(p)\otimes e_{pp} \\ \vspace{1\jot}
        X(s)&=\sum_{1\le p\le N}f(p;s)\otimes e_{pp}
          + \sum_{1\le p<q\le N}( g(p,q;s)\otimes e_{pq}+
            g(p,q;s)^\ast\otimes e_{qp}).
      \endalign
      $$
    Then $H$ and $X(s)$ are semicircular and
      $\{H,C,( X(s))_{s\in S}\}$ is $\ast$--free.
    Also, the von~Neumann algebra generated by $C$ and $H$ contains
      $1\otimes M_{N}(\Cpx)$ if the normal elements in $\nu_{1}$
      have disjoint spectra.
  \endproclaim

  \demo{Proof}
    By changing $A$ if necessary, we may assume that there exist
      a $\ast$--free family
      $\nu=\{U\}\cup\{V_{p}\mid2\le p\le N\}$ of
      Haar unitaries such that
      $\{\nu,\nu_{2},\nu_{3},\nu_{4},\nu_{5},\nu_{6}\}$
      is $\ast$--free,
      and bounded measurable functions $h_{1},\ldots,h_{N}$ on the unit circle
      such that $a(1)=h_{1}(U)$ and $a(p)=h_{p}(V_{p}^\ast UV_{p})$
      for $2\le p\le N$.
    Let
      $$ W=1\otimes e_{11}+\sum_{2\le p\le N}V_{p}\otimes e_{pp}. $$
    We shall show that $WHW^\ast$ is semicircular and that
      $\{WHW^\ast,WCW^\ast\}\cup\{WX(s)W^\ast\mid s\in S\}$
      is $\ast$--free.
    Note that
      $$ \align
        WCW^\ast&=\sum_{1\le p\le N}h_{p}(U)\otimes e_{pp} \\ \vspace{1.2\jot}
        WHW^\ast&= {\topaligned f(1)\otimes e_{11}
          &+\sum_{2\le p\le N}V_{p}f(p)V_{p}^\ast\otimes e_{pp}
            +\sum_{2\le q\le N}( b(q)V_{q}^\ast\otimes e_{1q}
            +V_{q}b(q)\otimes e_{q1}) \\ \vspace{0.7\jot}
          &+\sum_{2\le p<q\le N}( V_{p}g(p,q)V_{q}^\ast\otimes e_{pq}
            +V_{q}g(p,q)^\ast V_{p}^\ast\otimes e_{qp}) \endaligned} \\
            \vspace{1.2\jot}
        WX(s)W^\ast&={\topaligned f(1;s)\otimes e_{11}
          &+\sum_{2\le p\le N}V_{p}f(p;s)V_{p}^\ast\otimes e_{pp} \\
            \vspace{0.7\jot}
          &+\sum_{2\le q\le N}( g(1,q;s)V_{q}^\ast\otimes e_{1q}
            +V_{q}g(1,q;s)^\ast\otimes e_{q1}) \\ \vspace{0.7\jot}
          &+\sum_{2\le p<q\le N}( V_{p}g(p,q;s)V_{q}^\ast\otimes e_{pq}
            +V_{q}g(p,q;s)^\ast V_{p}^\ast\otimes e_{qp}). \endaligned}
      \endalign
      $$

    By Proposition~3.1, to prove
      semicircularity of $WHW^\ast$ and also $\ast$--freeness
      of
      $\bigl(WHW^\ast,WCW^\ast,\allowmathbreak(WX(s)W^\ast)_{s\in S}\bigr)$,
      it suffices to show that in $(A,\phi)$,
      $$
        \nu_{1}=\{f(1)\}\cup\{V_{p}f(p)V_{p}^\ast\mid2\le p\le N\}
        \cup\{f(1;s)\mid s\in S\}
        \cup\{V_{p}f(p;s)V_{p}^\ast\mid2\le p\le N,s\in S\}
      $$
      is a free family of semicircular elements,
      $$ \align
          \nu_{2}=
          &\{b(q)V_{q}^\ast\mid2\le q\le N\}
            \cup\{V_{p}g(p,q)V_{q}^\ast\mid2\le p<q\le N\}\cup
            \\ \vspace{0.5\jot}
          \cup&\{g(1,q;s)V_{q}^\ast\mid2\le q\le N,s\in S\}
            \cup\{V_{p}g(p,q;s)V_{q}^\ast\mid2\le p<q\le N,s\in S\}
        \endalign $$
      is a $\ast$--free family of circular elements and that
      $\{\nu_{1},\nu_{2},\{U\}\}$ is $\ast$--free.
    But (proceeding as in Voiculescu's proof of Theorem~3.3
      in~\cite{10}),
      we have polar
      decompositions
      $$ \align
        V_pg(p,q;s)V_q^\ast &= (V_pV(p,q;s)V_q^\ast)\,(V_q|g(p,q;s)|V_q^\ast)
\\
        V_pg(p,q)V_q^\ast &= (V_pV(p,q)V_q^\ast)\,(V_q|g(p,q)|V_q^\ast) \\
        \text{and }
          g(1,q;s)V_q^\ast &= (V(1,q;s)V_q^\ast)\,(V_q|g(1,q;s)|V_q^\ast),
      \endalign $$
      (where $g(p,q;s)=V(p,q;s)|g(p,q;s)|$
      and $g(p,q)=V(p,q)|g(p,q)|$ are polar decompositions).
    Now it is a question of showing that
      $$ \align
        \nu_1&\cup\{b(q)\mid2\le q\le N\} \\ \vspace{0.5\jot}
        &\cup\{V_q|g(p,q)|V_q^\ast\mid2\le p<q\le N\}\cup \\ \vspace{0.5\jot}
        &\cup\{V_pV(p,q)V_q^\ast\mid2\le p<q\le N\} \\ \vspace{0.5\jot}
        &\cup\{V_q\mid2\le q\le N\}\cup \\ \vspace{0.5\jot}
        &\cup\{V_q|g(p.q;s)|V_q^\ast\mid1\le p<q\le N,\,s\in S\}\cup
          \\ \vspace{0.5\jot}
        &\cup\{V(1,q;s)V_q^\ast\mid2\le q\le N\}\cup \\ \vspace{0.5\jot}
        &\cup\{V_pV(p,q;s)V_q^\ast\mid2\le p<q\le N,\,s\in S\}\cup
          \\ \vspace{0.5\jot}
        &\cup\{U\}
      \endalign $$
      is $\ast$--free.
    This is accomplished by using the functional calculus to
      replace $f(p)$, $f(p;s)$, $b(q)$, $|g(p,q)|$ and $|g(p,q;s)|$
      with equivalent Haar--unitaries (i.e\. Haar--unitaries which
      generate the same von Neumann algebra).
    Keeping Theorem~1.2 in mind, this reduces the question to one of whether
      certain elements of a free group form a set of free generators
      of a subgroup,
      which has in turn
      a straightforward solution.

    Now supposing that the normal elements in $\nu_{1}$ have disjoint
      spectra, we see that $1\otimes e_{pp}$ are spectral projections of
      $C$.
    But then $b(q)\otimes e_{1q}$ is in the $\ast$--algebra
      generated by $C$ and $H$,
      and the unitary part of its polar decomposition is $1\otimes e_{1q}$.
  \QED

  One can now prove results similar to the following.

  \proclaim{Corollary 3.6}
    $$
      L(\freeF_{N^{2}K})\otimes M_{N}(\Cpx)\cong L(\freeF_{K})
      \ast L^{\infty}(\{0\}\cup(1/N,1)),
    \tag36
    $$
    where $L^{\infty}(\{0\}\cup(1/N,1))$ is $L^{\infty}(\gamma)$
    equipped with trace $\int\cdot\dif\gamma$, and
    $\gamma$ is the probability measure having mass $1/N$ on $\{0\}$ and
    Lebesque measure on the open interval.
  \endproclaim

  \demo{Proof}
    Let $S$ be a set of cardinality $K-1$.
    Then let $L(\freeF_{N^{2}K})\subseteq A$ be generated by a free family of
      Haar unitaries, $\nu_{1}=\{\tilde{u_{2}},\cdots,\tilde{u_{N}}\}$
      together with the families $\nu_{2},\cdots,\nu_{6}$ of Proposition~3.5.
    Let
      $$ C=1\otimes e_{11}+\sum_{2\le p\le N}p\tilde{U_{p}}\otimes e_{pp}. $$
    Then $C$ generates the $L^{\infty}$ space while
      $H$ and $\{X(s)\mid s\in S\}$ generate $L(\freeF_{K})$.
    But since, by the proposition, $C$ and $H$ together generate the matrix
      units, we see that $C$, $H$, and $\{X(s)\mid s\in S\}$ generate
      all of the left--hand side of~(36).
  \QED

  Note that in exactly the same way, one obtains
    $$ L(\freeF_{N^2K-N+1+j})\otimes M_N(\Cpx)
      \cong L(\freeF_K)\ast L^\infty(\{0\}\cup\cdots\cup\{N-1-j\}\cup
      (\tfrac{N-j}N,1)) \tag37 $$
    for $0\le j\le N$, where again each atom has measure $1/N$
    and the interval has Lebesque measure.
  In~\cite{10}, Voiculescu did
     the cases $j=0$ and $j=N$.

\noindent{\bf\S4\. Free products of $L(\freeF_K)$
  with the hyperfinite II$_1$--factor}

  Let $R$ denote the hyperfinite II$_{1}$--factor.
  Since $L(\freeF_K)\ast R$ is an inductive limit
    of $L(\freeF_K)\ast M_{2^{k}}(\Cpx)
    \cong L(\freeF_K)\ast L^{\infty}(\{0\}\cup(1/2^k,1))$
    and $L^{\infty}(\{0\}\cup(1/2^k,1))$ ``approaches'' $L(\Integers)$
    as $k\rightarrow\infty$,
    we could hope that
    $L(\freeF_K)\ast R\cong L(\freeF_{K+1})$.
  This is indeed the case, however, as the natural generators for the
    $L^{\infty}(\{0\}\cup(1/2^{k},1))$'s don't
    map nicely to each other under
    the inductive limit homomorphisms,
    it seemed intractable to try a proof using inductive limits
    and thus the following proof, though it is, in the end,
    an inductive argument, consists of
    simply exhibiting the requisite free generators.

  \proclaim{Theorem 4.1}
    For $1\le K\le\infty$,
      $$ L(\freeF_K)\ast R\cong L(\freeF_{K+1}), $$
      where $R$ is the hyperfinite II$_{1}$--factor.
  \endproclaim

  \demo{Proof}
    Since $\ast$ is associative, it suffices to
      show $L(\Integers)\ast R\cong L(\freeF_{2})$.
    Let $X$ be a semicircular random variable generating $L(\Integers)$.
    In the hyperfinite II$_{1}$--factor, $R$, let
      $\Omega_{k}$ be the $2^{k}\times2^{k}$ system of matrix units,
      $$ \Omega_{k}=\{ e(p_{1}\cdots p_{k},q_{1}\cdots q_{k};2^{k})
        \mid p_{i},q_{i}\in\{0,1\}\}, $$
      for $k\ge1$,
      where $p_{1}\cdots p_{k}$ and $q_{1}\cdots q_{k}$ are binary
      expansions, such that
      $$ e(p_{1}\cdots p_{k},q_{1}\cdots q_{k};2^{k})
        =e(p_{1}\cdots p_{k}0,q_{1}\cdots q_{k}0;2^{k+1})
        +e(p_{1}\cdots p_{k}1,q_{1}\cdots q_{k}1;2^{k+1}) $$
      and $\bigcup_{k\ge1}\Omega_{k}$ generates $R$.
    Let $\MvN=L(\Integers)\ast R$ with trace $\phi$.

    Let $\Ac_{k}$ denote the algebra
      $e(\OdO,\OdO;2^{k})\MvN e(\OdO,\OdO;2^{k}),$
      containing
      $$ \align
        F(p_{1}\cdots p_{k};k)
          &=2^{k/2}e(\OdO,p_{1}\cdots p_{k};2^{k})
          \,X\,e(p_{1}\cdots p_{k},\OdO;2^{k}) \\ \vspace{0.5\jot}
        G(p_{1}\cdots p_{k},q_{1}\cdots q_{k};k)
          &=2^{k/2}e(\OdO,p_{1}\cdots p_{k};2^{k})
          \,X\,e(q_{1}\cdots q_{k},\OdO;2^{k}).
      \endalign $$
    Then by Proposition~3.2,
      $$ \nu_{1}=\{ F(p_{1}\cdots p_{k};k)\mid p_{i}\in\{0,1\}\} $$
      is a free semicircular family and
      $$ \nu_{2}=\{ G(p_{1},\cdots p_{k},q_{1}\cdots q_{k};k)
        \mid p_{i},q_{i}\in\{0,1\}, p_{1}\cdots p_{k}>q_{1}\cdots q_{k}\} $$
      a $\ast$--free circular family in $(\Ac_{k},2^k\phi\restrict_{\Ac_k})$
      such that $\{\nu_{1},\nu_{2}\}$ is
      $\ast$--free.
    If $a\in \Ac_{k}$ let
      $a\otimes e(p_{1}\cdots p_{k},q_{1}\cdots q_{k};2^{k})$ denote
      $e(p_{1}\cdots p_{k},\OdO;2^{k})\,a\,e(\OdO,q_{1}\cdots q_{k};2^{k})
      \in\MvN$.

    Let $G(\OdO1,\OdO;k)=V(k)B(k)$ be the polar decomposition,
      and let
      $$ \align
        H&={\topaligned \sum_{k=1}^{\infty}2^{-k/2}
          \bigl(&F(\OdO1;k)\otimes e(\OdO1,\OdO1;2^{k}) \\
            &+B(k)\otimes( e(\OdO,\OdO1;2^{k})
            +e(\OdO1,\OdO;2^{k}))\bigr) \endaligned} \\ \vspace{1\jot}
        W&=\sum_{k=1}^{\infty}\frac{1}{k}V(k)\otimes e(\OdO1,\OdO1;2^{k}),
      \endalign $$
      (see Figure~5).

    It is not hard to show that
      $H$ and $W$ together generate $L(\Integers)\ast R$.
    Indeed, from spectral projections of $W$, one has many of the
      diagonal matrix units, enough to extract from
      $H$ the operators, $B(k)\otimes e(\OdO,\OdO1;2^k)$ for $1\le k<\infty$.
    Taking the polar decompositions of these, we get as polar parts
      $1\otimes e(\OdO,\OdO1;2^k)$.
    The ensemble of these, together with the aforementioned spectral
      projections of $W$, generate all the matrix units, hence $R$.
    Now that we have the matrix units, we can pull apart
      $W$ and $H$, and piece the parts together to obtain $X$
      (as a norm limit).

    Now we shall show that $H$ is semicircular and $\{W,H\}$ is $\ast$--free.
    This will prove the theorem because the spectral measure of $W$ has no
      atoms, so each of $W$ and $H$ will generate a copy of $L(\Integers)$.
    For $n\ge1$, let
      $$ \align
        H_{n}&=\left(\sum_{k=1}^{n}2^{-k/2}
          {\topaligned(\,&F(\OdO1;k)\otimes e(\OdO1,\OdO1;2^{k}) \\
            &+B(k)\otimes( e(\OdO,\OdO1;2^{k})
            +e(\OdO1,\OdO;2^{k}))\,)\,\endaligned}\right) \\
          &+2^{-n/2}F(\OdO;n)\otimes e(\OdO,\OdO;2^{n}) \\ \vspace{1\jot}
        W_{n}&=\left(\sum_{k=1}^{n}\frac{1}{k}V(k)\otimes e(\OdO1,\OdO1;2^{k})
          \right)+\left(\tfrac{1}{n+1}\right)1\otimes e(\OdO,\OdO;2^{n}).
      \endalign $$
    Since $F(\cdots;k)$ has the same norm for all $k$,
      and similarly for $G(\cdots;k)$,
      clearly $W_{n}$ and $H_{n}$ converge in norm to $W$ and $H$.
    Hence it suffices to show that for each $n$, $H_{n}$ is semicircular
      and $\{W_{n},H_{n}\}$ is $\ast$--free.
    This we do by induction on $n$.
    In the case $n=1$ we have
      $$ H_{1}= 2^{-1/2}
        \left(
        \matrix
          F(0;1) & B(1) \\
          B(1) & F(1;1)
        \endmatrix
        \right),
        \quad W_{1}=
        \left(
        \matrix
          1/2 & 0 \\
          0 & V(1)
        \endmatrix
        \right). $$
    Then by Propositions~3.5 and~3.2 and Theorem~1.2, $H_{1}$ is
      semicircular and $\{W_{1},H_{1}\}$ is $\ast$--free.

    For the inductive step, suppose the proposition is true for $n-1$.
    Let
      $$ \align
        K_{n}&=2^{1/2}(H_{n}-F(1;1)\otimes e(1,1;2)
          -B(1)\otimes( e(1,0;2)+e(0,1;2))\,) \\ \vspace{0.5\jot}
        Y_{n}&=W_{n}-V(1)\otimes e(1,1;2).
      \endalign $$
    Then by inductive hypothesis applied to $\Ac_{1}$, we have that
      in $(\Ac_{1},2\phi\restrict_{\Ac_1})$,
      $K_{n}$ is semicircular and $\{K_{n},Y_{n}\}$ is $\ast$--free.
    Note that
      $$ H_{n}= 2^{-1/2}
        \left(
        \matrix
          K_{n} & B(1) \\
          B(1) & F(1;1)
        \endmatrix
        \right),
        \quad W_{n}=
        \left(
        \matrix
          Y_{n} & 0 \\
          0 & V(1)
        \endmatrix
        \right). $$

    In $(\Ac_{1},2\phi)$,
      let $S=\{ e(0p_{1}\cdots p_{k-1},0q_{1}\cdots q_{k-1};2^{k})
      \mid p_{i},q_{i}\in\{0,1\},\,2\le k\le n\}$.
    Then by Proposition~3.3,
      $\{S,F(0;1),G(1,0;1),F(1;1)\}$ is $\ast$--free.
    But $K_{n}$ and $Y_{n}$ are generated by $S$ and $F(0;1)$, so
      using Proposition~1.1 and Theorem~1.2,
      we obtain that $\{K_{n},Y_{n},V(1),B(1),F(1;1)\}$ is $\ast$--free.
    Thus by Proposition~3.5,
      $H_{n}$ is semicircular and $\{W_{n},H_{n}\}$ is $\ast$--free.
  \QED

  \proclaim{Corollary 4.2}
    For $1\le n<\infty$, let $(G_i)_{i=1}^n$ be a family of discrete
      groups such that for each $i$, $G_i$ is either an amenable i.c.c\.
      group or a copy of $\Integers$, and at least one of them is the latter.
    Then
      $$ L(\freeprod_{i=1}^nG_i)\cong L(\freeF_n). \tag38 $$

    Moreover, suppose now that $(G_i)_{i=1}^\infty$ is a family of
      amenable i.c.c\. discrete groups.
    Then
      $$ L(\freeF_\infty\ast(\freeprod_{i=1}^\infty G_i))
        \cong L(\freeF_\infty). \tag39 $$
  \endproclaim

  \demo{Proof}
    Let us first recall Connes' result~\cite{3}
      that whenever $G$ is an amenable i.c.c\.
      discrete group, $L(G)$ is a copy of $R$, the hyperfinite
      II$_1$--factor.
    Now~(38) follows immediately.
    For~(39), note that free products may be taken in any order, so
      $$ L(\freeF_\infty\ast(\freeprod_{i=1}^\infty G_i))
        \cong L(\freeprod_{i=1}^\infty(\Integers\freeprod G_i))
        \cong \freeprod_{i=1}^\infty L(\Integers\freeprod G_i)
        \cong \freeprod_{i=1}^\infty L(\freeF_2)
        \cong L(\freeF_\infty). $$
  \QED

  \noindent{\bf Acknowledgements.}
  I would like to thank Dan Voiculescu, my advisor,  for many
    helpful discussions, and also the people at the Centre
    de Recherches Math\'{e}matiques at the University of Montreal,
    where part of this work was done, for providing a friendly
    and stimulating atmosphere during their special semester
    in operator algebras.

\vskip4ex
\noindent{\bf Appendix}

  As F\. R\u{a}dulescu observed in~\cite{17},
    there are proofs of Gaussian cases of Theorem~2.1 which are very similar
    to Voiculescu's proofs for his random matrix results in~\cite{11}.
  Consider for example the two cases having the following
    additional hypotheses for Theorem~2.1.
  \roster
  \item"(I)" Real Gaussian:  each $a(i,j;s,n)$ is a real Gaussian
    random variable and $\{a(i,j;s,n)\mid s\in S,\,1\le i\le j\le n\}$
    is an independent family.
  \item"(II)" Complex Gaussian: each $a(k,k;s,n)$ is a real Gaussian
    random variable, each $\Rlpt a(i,j;s,n)$ and $\Imaginary a(i,j;s,n)$
    for $i\ne j$ are real Gaussian random variables having second
    moments equal to $1/2n$ and
    $$ \{\Rlpt a(i,j;s,n)\mid s\in S,\,1\le i\le j\le n\}
      \cup\{\Imaginary a(i,j;s,n)\mid s\in S,\,1\le i<j\le n\} $$
    is an independent family.
  \endroster

  We remark that the complex Gaussian case alone
    suffices to prove Proposition~3.2.

  In order to prove~2.1 for these two cases, one proceeds exactly as in the
    proof of Theorem~2.2 of~\cite{11}, which utilizes Theorem~2.1 of~\cite{11}
    and the fact that a sum of Gaussian random variables is Gaussian to reduce
    the proof to checking certain second order freeness conditions.
  These are in turn checked in our case of having block diagonal constant
    matrices by using counting arguments only slightly more complicated than
    those in~\cite{11}.

  \vskip7ex
  An earlier version of this paper was distributed as
    preprint~CRM--1750.

\end{document}

__________________________________________________________________________

%
%

\documentstyle{article}

\newcommand{\vtilde}{\stackrel{v}{\sim}}         

\begin{document}

    \begin{figure}[hbt]
    \centering
      \begin{picture}(250,29)(5,-12)
        \put(0,0){\line(1,0){120}}
        \put(168,0){\line(1,0){72}}
        \put(0,0){\circle*{4}}
        \put(48,0){\circle*{4}}
        \put(96,0){\circle*{4}}
        \put(192,0){\circle*{4}}
        \put(240,0){\circle*{4}}
        \put(132,0){\circle*{2}}
        \put(144,0){\circle*{2}}
        \put(156,0){\circle*{2}}
        \put(-6,6){$1$}
        \put(47,7){$2$}
        \put(96,7){$3$}
        \put(187,7){$m$}
        \put(235,6){$m+1$}
        \put(14,-11){$\scriptstyle a(i_{1},j_{1})$}
        \put(62,-11){$\scriptstyle a(i_{2},j_{2})$}
        \put(198,-11){$\scriptstyle a(i_{m},j_{m})$}
      \end{picture}
    \vspace{2ex}
    \caption{The straight line graph of length $m$.}
    \label{figure:straight_line_graph}
    \end{figure}

    \begin{figure}[hbt]
    \centering
      \begin{picture}(197,175)(5,-58)
        \put(0,0){\line(1,0){168}}
        \put(96,0){\line(0,-1){48}}
        \put(0,0){\circle*{4}}
        \put(48,0){\circle*{4}}
        \put(96,0){\circle*{4}}
        \put(144,0){\circle*{4}}
        \put(96,-48){\circle*{4}}
        \put(180,0){\circle*{2}}
        \put(192,0){\circle*{2}}
        \put(204,0){\circle*{2}}
        \put(-6,6){$1$}
        \put(46,7){$2$}
        \put(86,7){$3\vtilde5$}
        \put(142,7){$6$}
        \put(99,-58){$4$}
      \end{picture}
    \vspace{2ex}
    \caption{A click.}
    \label{figure:click}
    \end{figure}

    \begin{figure}
    \centering
      \begin{picture}(120,180)

        \put(0,0){\line(0,1){2.5}}
        \multiput(0,7.5)(0,15){7}{\line(0,1){5}}
        \put(0,117.5){\line(0,1){2.5}}

        \put(40,0){\line(0,1){2.5}}
        \multiput(40,7.5)(0,15){7}{\line(0,1){5}}
        \put(40,117.5){\line(0,1){2.5}}

        \put(80,0){\line(0,1){2.5}}
        \multiput(80,7.5)(0,15){7}{\line(0,1){5}}
        \put(80,117.5){\line(0,1){2.5}}

        \put(120,0){\line(0,1){2.5}}
        \multiput(120,7.5)(0,15){7}{\line(0,1){5}}
        \put(120,117.5){\line(0,1){2.5}}

        \put(0,0){\line(1,0){2.5}}
        \multiput(15,0)(15,0){7}{\line(1,0){5}}
        \put(117.5,0){\line(1,0){2.5}}

        \put(0,40){\line(1,0){2.5}}
        \multiput(15,40)(15,0){7}{\line(1,0){5}}
        \put(117.5,40){\line(1,0){2.5}}

        \put(0,80){\line(1,0){2.5}}
        \multiput(15,80)(15,0){7}{\line(1,0){5}}
        \put(117.5,80){\line(1,0){2.5}}

        \put(0,120){\line(1,0){5}}
        \multiput(15,120)(15,0){7}{\line(1,0){5}}
        \put(117.5,120){\line(1,0){2.5}}

        \put(0,0){\line(1,1){120}}

        \put(-35,-10){$(0,0)$}
        \put(126,122){$(3,3)$}

      \thicklines

        \put(0,0){\line(0,1){40}}
        \put(0,40){\line(1,0){40}}
        \put(40,40){\line(0,1){80}}
        \put(40,120){\line(1,0){80}}


      \end{picture}
    \vspace{2ex}
    \caption{The path for $(-1,1,-1,-1,1,1)$}
    \end{figure}

    \begin{figure}[hbt]
    \centering
      \begin{picture}(160,210)(-30,-48)
        \put(0,0){\line(1,1){33.94}}
        \put(33.94,33.94){\line(1,0){48}}
        \put(81.94,33.94){\line(1,-1){33.94}}
        \put(115.88,0){\line(-1,-1){33.94}}
        \put(33.94,-33.94){\line(-1,1){33.94}}
        \put(0,0){\circle*{4}}
        \put(33.94,33.94){\circle*{4}}
        \put(81.94,33.94){\circle*{4}}
        \put(115.88,0){\circle*{4}}
        \put(81.94,-33.94){\circle*{4}}
        \put(33.94,-33.94){\circle*{4}}
        \put(69.94,-33.94){\circle*{2}}
        \put(57.94,-33.94){\circle*{2}}
        \put(45.94,-33.94){\circle*{2}}
        \put(-12,0){$m$}
        \put(25.94,37.94){$1$}
        \put(84.94,37.94){$2$}
        \put(118.88,0){$3$}
        \put(84.94,-48.94){$4$}
        \put(20.94,-48.94){$m-1$}
        \put(118.88,-48.94){.}
        \put(42,49){$a(i_{1},j_{1})$}
        \put(110,29){$a(i_{2},j_{2})$}
        \put(-30,29){$a(i_{m},j_{m})$}
        \put(110,-29){etc.}
      \end{picture}
    \vspace{2ex}
    \caption{The m-gon graph.}
    \label{figure:m-gon_graph}
    \end{figure}

    \begin{figure}[p]
    \setlength{\unitlength}{4in}
    \begin{picture}(1,1)
      \put(0,0.5){\line(1,0){1}}
      \put(0,0.75){\line(1,0){0.5}}
      \put(0,0.875){\line(1,0){0.25}}
      \put(0.5,0){\line(0,1){1}}
      \put(0.25,0.5){\line(0,1){0.5}}
      \put(0.125,0.75){\line(0,1){0.25}}
      \put(0.175,0.245){\LARGE $B(1)$}
      \put(0.675,0.745){\LARGE $B(1)$}
      \put(0.700,0.245){\LARGE $F(1;1)$}
      \put(0.07,0.625){\large$B(2)$}
      \put(0.32,0.875){\large$B(2)$}
      \put(0.29,0.625){\large$F(01;2)$}
      \put(0.025,0.812){$B(3)$}
      \put(0.15,0.932){$B(3)$}
      \put(0.1425,0.812){$F(\cdots)$}
      \put(0.03125,0.96875){\Huge$.$}
      \put(0.0625,0.9375){\Huge$.$}
      \put(0.09375,0.90625){\Huge$.$}
    \end{picture}
    \vspace*{1.5ex}
    \caption{A picture of $H$ in the proof of Theorem 4.1}
    \label{figure:H}
    \end{figure}


\Refs

  \ref \no 1 \by L\. Arnold \paper Deterministic version
    of Wigner's semicircle law for the distribution of matrix eigenvalues
    \jour Lin\. Algebra Appl\. \vol 13 \yr 1976 \pages 185--199 \endref

  \ref \no 2 \by W.M\. Ching \paper Free products of von Neumann algebras
    \jour Trans\. Amer\. Math\. Soc\. \vol 178 \yr 1973 \pages 147--163 \endref

  \ref \no 3 \by A\. Connes \paper Classification of injective factors
    \jour Ann\. of Math\. \vol 104 \yr 1976 \pages 73--115 \endref

  \ref \no 4 \by M.L\. Mehta \book Random Matrices \bookinfo Revised and
    enlarged 2nd ed\. \publ Academic Press \yr 1991 \endref

  \ref \no 5 \by F.J\. Murray and J. von Neumann
    \paper Rings of operators\.~IV
    \jour Ann\. of Math\. \vol 44 \yr 1943 \pages 716--808 \endref

  \ref \no 6 \by F\. R\u{a}dulescu \paper
    The fundamental group of the von Neumann algebra of a free group with
    infinitely many generators is $\Reals_{+}\backslash\{0\}$
    \jour J\. Amer\. Math\. Soc\. \toappear \endref

  \ref \no 7 \by J.W\. Silverstein \paper Eigenvalues
    and Eigenvectors of large dimensional sample covariance matrices
    \jour Contemporary Math\. \vol 50 \yr 1986 \pages 153--159 \endref

  \ref \no 8 \manyby D\.  Voiculescu \paper Symmetries of some
    reduced free product C$^{\ast}$--algebras \inbook Operator Algebras
    and Their Connections with Topology and Ergodic Theory  \publ Lecture
    Notes in Mathematics, Volume~1132, Springer--Verlag \yr 1985
    \pages 556--588 \endref

  \ref \no 9 \bysame \paper Noncommutative random variables and spectral
    problems in free product C$^*$--algebras \jour Rocky Mt\. J. Math\.
    \vol 20 \yr 1990 \pages 263-283 \endref

  \ref \no 10 \bysame \paper Circular and semicircular systems
    and free product factors \inbook Operator Algebras, Unitary
Representations,
    Enveloping Algebras, and Invariant Theory \bookinfo Progress in
    Mathematics, Volume~92, \publ Birkh\"{a}user
    \publaddr Boston \yr 1990 \endref

  \ref \no 11 \bysame \paper Limit laws for Random
    matrices and free products \jour Invent\. Math\. \vol 104
    \yr 1991 \pages 201-220 \endref

  \ref \no 12 \bysame \paper Free non-commutative random
    variables, random matrices and the II$_1$--factors
    of free groups \inbook Quantum Probability and Related
    Topics~VI \bookinfo L\. Accardi, ed\. \publ World Scientific
    \publaddr Singapore \yr 1991 \pages 473--487 \endref

  \ref \no 13 \by K.W\. Wachter \paper The strong limits
    of random matrix spectra for sample matrices of
    independent elements \jour Ann\. Prob\. \vol 6 \yr 1978
    \pages 1--18 \endref

  \ref \no 14 \manyby E.P\. Wigner \paper Characteristic
    vectors of bordered matrices with infinite dimensions
    \jour Ann\. of Math\. \vol 62 \yr 1955 \pages 548--564 \endref

  \ref \no 15 \bysame \paper On the distribution
    of the roots of certain symmetric matrices
    \jour Ann\. of Math\. \vol 67 \yr 1958 \pages 325--327 \endref

  \ref \no 16 \by A.M\. Yaglom and I.M\. Yaglom \book Challenging problems
    with elementary solutions \bookinfo Translated by James, McCawley
    Jr., Revised and edited by Basil Gordon \publaddr San Francisco
    \publ Holden--Day \yr 1964 \endref

  \ref \no 17 \by F\. R\u{a}dulescu \paper Random matrices,
    amalgamated free products
    and subfactors of the von Neumann algebra of a free group
    \paperinfo preprint, I.H.E.S\., December, 1991 \endref

\endRefs
\end{document}